\documentclass[reqno]{amsart}

\usepackage{amssymb,amsfonts,amsmath,mathdots,amsthm}
\usepackage{float}

\usepackage{tikz-cd}
\usepackage{lmodern}
\usepackage{mathrsfs}
\usepackage[sc]{mathpazo}
\DeclareMathAlphabet{\mathfrak}{U}{euf}{m}{n}
\usepackage{mathtools}  %xhook, \ev
\usepackage{calligra}
\usepackage{xcolor}
\usepackage[hidelinks]{hyperref}
\usepackage{accents}
\usepackage{subcaption}
\usepackage{enumitem}

\usepackage{hyperref}

\hypersetup{colorlinks=true,
	linkcolor=magenta,
	citecolor=magenta}

\setlist[itemize]{leftmargin=*}

\newcommand\zero{0}

\usepackage{relsize}
\usepackage{etoolbox,fancyhdr}
\usepackage{indentfirst}  

\usetikzlibrary{arrows,arrows.meta,shapes.misc,decorations.markings,positioning,patterns}

\mathtoolsset{showonlyrefs} %Only tag reffered equations

\DeclarePairedDelimiter\ev{\langle}{\rangle}%  No langle rangle

\DeclareMathOperator{\Z}{\mathbb{Z}}
\DeclareMathOperator{\R}{\mathbb{R}}

\DeclareMathOperator{\Tr}{Tr}
\DeclareMathOperator{\Res}{Res}

\DeclareMathOperator{\Tor}{Tor}

\DeclareMathOperator{\colimit}{\textup{colim}}

\theoremstyle{plain}

\newtheorem{prop}{Proposition}[section]

\newtheorem{conj}[prop]{Conjecture}
\makeatother
\title{The $RO(C_4)$ integral homology of a point}
\author{Nick Georgakopoulos}

\begin{document}

\maketitle{}

\begin{abstract}We compute the $RO(C_4)$ integral homology of a point with complete information as a Green functor, and we show that it is generated, in a slightly generalized sense, by the Euler and orientation classes of the irreducible real $C_4$-representations. We have devised a computer program that automates these computations for groups $G=C_{p^n}$ and we have used it to verify our results for $G=C_4$ in a finite range.	
\end{abstract}

\setcounter{tocdepth}{1}
\tableofcontents
\section{Introduction}\label{Intro}

The computation of the $RO(G)$-graded (Bredon) homology of a point has been historically a very difficult problem. Stong and Lewis completely determined it for $G=C_p$, the cyclic group of prime order $p$, using coefficients in the Burnside Green functor (c.f. \cite{Lew88}). More recently, there has been renewed interested in such computations, but now for non prime-order groups $G$ and using constant $\underline{\Z}$ coefficients. The reason for this resurgence is rooted in the seminal work \cite{HHR16}, where $RO(C_8)$-graded homology is used to solve the Kervaire invariant problem in all but one dimension.

The very first step in \cite{HHR16} is obtaining partial information about the $RO(C_8)$ integral homology of a point, 
\begin{equation}
\underline{\pi}_{\bigstar}^{C_8}(H\underline \Z)=\underline{H}_{\bigstar}^{C_8}(S^0;\underline \Z).
\end{equation} 
The vanishing of these homologies for certain virtual representations $\bigstar$ is used to prove the Gap Theorem. In \cite{HHR17}, they perform a similar computation in the more tractable $C_4$ case, obtaining more complete information for the slice spectral sequence they consider. However, because they are interested only in the integer-graded part of that spectral sequence, they only compute $\underline{H}_{\bigstar}^{C_4}(S^0;\underline{\Z})$ for $\bigstar=k+r\rho$ where $k,r\in \Z$ and $\rho$ is the regular $C_4$-representation.

The computational aspects of equivariant homotopy theory have been much
less explored than their nonequivariant counterparts, and that’s partly
due to the greater complexity of the algebra involved.  As we illustrate in this paper, this complexity is reflected in the answers retrieved by the calculations and is not just a technicality or a limitation of our approach. The methods we
use to obtain the calculations are both theoretical and computer based, and can be applied to a greater range of equivariant computations.

In this paper, we compute $\underline{H}_{\bigstar}^{C_4}(S^0;\underline{\Z})$ for all possible $\bigstar$, as a Green functor. This means that in addition to all the groups $H_k^H(S^V;\underline{\Z})$ for $V$ a real virtual $C_4$ representation and $H$ a subgroup of $C_4$, together with their restrictions and transfers, we also compute the multiplicative structure which comes from the fact that $H\underline{\Z}$ is a $C_4$-ring spectrum.

The group $C_4$ has two non-trivial irreducible real representations, the $1$-dimen\-sional sign representation $\sigma$ and the $2$-dimensional representation $\lambda$  (rotation by $\pi/2$ degrees). Therefore, we are in effect computing the Mackey functors $\underline{H}_k^{C_4}(S^{\pm n\sigma\pm m\lambda};\underline{\Z})$ for $k,n,m\in \Z$ and $n,m\ge 0$. When both signs in $S^{\pm n\sigma\pm m\lambda}$ are positive, we have an explicit and simple equivariant cellular decomposition for the space $S^{n\sigma+m\lambda}$ and we can compute the homology using the cellular chain complex $\underline C_*S^{n\sigma+m\lambda}$. When both signs are negative, we can appeal to Spanier Whitehead Duality: 
\begin{equation}
\underline H_k^{C_4}(S^{-n\sigma-m\lambda})=\underline H^{-k}_{C_4}(S^{n\sigma+m\lambda})
\end{equation}
and this is the cohomology of the cochain complex $\underline C^{-*}(S^{n\sigma+m\lambda})$ dual to the chains $\underline C_{-*}(S^{n\sigma+m\lambda})$ we had before.

The more difficult part of the computation is when we have opposite signs, such as $\underline H_k^{C_4}(S^{n\sigma-m\lambda};\underline{\Z})$. In this case, we could in principle work with the box product of chain complexes 
\begin{equation}
\underline C_*S^{n\sigma}\boxtimes \underline C^{-*}S^{m\lambda}
\end{equation}
but these complexes get intractably large for calculations by hand as $n,m$ get large. In place of that, we instead make use of three algebraic spectral sequences associated to these complexes: Two Atiyah-Hirzebruch spectral sequences and a Kunneth spectral sequence. Comparison of these three allows us to get the answer through fairly intuitive (if lengthy) arguments. A complication  is that everything needs to be performed on the Mackey functor level: for example, the $\Tor$ terms in the Kunneth spectral sequence are computed in the symmetric monoidal category of $\underline{\Z}$-modules. 

The main result of this paper is that $\underline{\pi}_{\bigstar}^{C_4}(H\underline \Z)$  is generated, in a generalized sense, by the Euler and orientation classes associated to $\sigma,\lambda$. These classes, under the operations of multiplication, division (see subsection \ref{Division} for the precise meaning of "division"), restriction and transfer, don't quite generate the entire $\underline{\pi}_{\bigstar}^{C_4}(H\underline{\Z})$, missing the generator of $H^{C_4}_{-3}(S^{-2\lambda})=\Z/4$. However, it turns out that this $\Z/4$ fits in a short exact sequence of abelian groups
\begin{equation}
0\to \Z/2\to \Z/4\to \Z/2\to 0
\end{equation}
where the $\Z/2$'s are obtainable using just the Euler and orientation classes. So if we include this group extension into our list of "operations", then the closure of the Euler and orientation classes under said operations is the entire $\underline{\pi}_{\bigstar}^{C_4}(H\underline{\Z})$. To be more precise, since we are interested in the homology as a Mackey functor, we shouldn't adjoin a group extension but rather the Mackey functor extensions that induce it.

At this point, we should mention earlier work by Zeng on this topic. \cite{Zeng} calculates the integer coefficient $RO(C_{p^2})$-graded homology of a point for all primes $p$, using the associated Tate-square diagram as opposed to the cellular chains approach we use here. His description for the multiplicative structure is in terms of the connecting homomorphism of certain cofiber sequences, while our description is solely in terms of the Euler and orientation classes. Modulo this difference, our results agree with his for the case $p=2$. We hope this write-up provides a more comprehensive analysis of this case, while also offering a detailed discussion of the subtleties involved in this computation, many of which are relevant in correctly interpreting the results of the computation (c.f. subsection \ref{Subtle}).

Another novelty in our work is the computerization of this computation, not just for $G=C_4$ but indeed for any $G=C_{p^n}$. We have devised a computer program that automatically produces the answer for both the additive and multiplicative structures of $\underline \pi_{\bigstar}^G(H\underline{\Z})$ or more generally $\underline \pi_{\bigstar}^G(H\underline{R})$ where $R$ is a user specified ring such as $\mathbb F_p$ or $\mathbb Q$. It can also compute the Massey products present in $\underline \pi_{\bigstar}^G(H\underline R)$ together with their indeterminacy. Of course the program can only work in a finite range, i.e. it can produce the answer for $S^{V}$ where the dimension of $V$ is bounded. 

Therefore, while the code can't completely replace the proof-based work, at the very minimum it's a powerful verification tool. For example, it was able to spot a few edge cases where mistakes were present in an early draft of this paper. The other advantage is scalability: determining $\underline \pi_{\bigstar}^{C_{p^n}}(H\underline{\Z})$ by hand for $n\ge 3$ is significantly more laborious than the $n=2$ case as there are more representations to contend with. But using the computer program, we can quickly and easily compute $\underline \pi_{\bigstar}^{C_{p^n}}(H\underline{\Z})$ in a large range and get a good grasp for what the answer should be. Indeed, the computational data our computer generated for groups $G=C_8,C_{16}$ and $C_{32}$ led us to formulate a conjecture for $\pi_{\bigstar}^{C_{2^n}}(H\underline{\Z})$ and all $n\ge 2$, that we describe in section \ref{Conjecture}. The spectral sequence computations in \cite{BC4S2} were also aided by computational results derived from our program. 

As for the organization of this paper, section \ref{Computer} offers a brief introduction of how our program works. The rest of the paper is completely independent of that.

Section \ref{MackeyFunctors} includes the 16 Mackey functors that appear in $\underline{H}_k^{C_4}(S^{\pm n\sigma\pm m\lambda})$ and their notation used throughout this paper. Remarkably we can get a non-cyclic answer for $H_k^{C_4}(S^{n\sigma-m\lambda})$ but only for even $n\ge 4$.

Section \ref{Generators} expounds on how the Euler and orientation classes generate $\underline{\pi}_{\bigstar}^{C_4}(H\underline{\Z})$ and how all relations can be effectively reduced to a single one, the "Gold Relation".

Section \ref{Conjecture} contains our general conjecture for the multiplicative structure of $\underline \pi_{\bigstar}^{C_{2^n}}(H\underline \Z)$ backed up by computational data for $n\le 5$.

Section \ref{Results} includes the complete determination of the Green functor $\underline{\pi}_{\bigstar}^{C_4}(H\underline{\Z})$ in the form of 8 readily usable tables.  

Section \ref{Framework} summarizes the theoretical framework of our computations.

The last five sections include the proofs of our results and make up the bulk of this paper. The final one is an Appendix devoted to proving that the Gold relation generates all other relations.

%We plan to investigate the $C_8$ case, and whether $\underline \pi_{\bigstar}^{C_8}(H\underline{\Z})$ is generated by the Euler and orientation classes (in the same generalized sense) in a sequel. If true, it would raise the question whether this continues to happen for $G=C_{2^n}$, $n\ge 4$. We also aim to investigate the case of $\underline{\mathbb F_2}$ coefficients and any interesting Massey products that arise, now aided by our computer program.

%The computational aspects of equivariant homotopy theory have been much less explored than their nonequivariant counterparts, and that's partly due to the greater complexity of the algebra involved. This complexity is not just a technicality, as it's reflected in the answers retrieved by the calculations. 
By computerizing the complicated algebra involved in these computations, we can drastically expand the known calculations while reducing the human work required. Currently, our computer program can do the $G=C_{p^n}$ case but in the future we expect to extend this to arbitrary finite abelian groups and certain nonabelian groups with known representation rings. The source code is publicly available \href{https://github.com/NickG-Math/Mackey}{here}, where the interested reader can not only inspect it, but also contribute to its improvement and expansion, which we highly encourage.

\subsection*{Acknowledgment} We want to thank Mingcong Zeng for carefully reading an earlier draft of this paper, helping us compare our computations with his, and for pointing out a subtlety in a certain relation of the multiplicative structure (c.f. end of subsection \ref{Subtle}). We would also like to thank Peter May for reading several earlier drafts of this paper.

\section{The computer program}\label{Computer}

The computations in this paper rely on filtering box products such as
\begin{equation}
C_*(S^{n\sigma-m\lambda})=C_*(S^{n\sigma})\boxtimes C_*(S^{-m\lambda})
\end{equation}
in different ways and comparing the resulting spectral sequences (here and always, $n,m$ are nonnegative integers). Ideally, we would  be working directly with that box product, but there are two major complications that prohibit this: Firstly, the box product of Mackey functors is not the level-wise tensor product. Instead, only the bottom level (corresponding to the orbit $G/e$) can be obtained as the tensor product, while all the higher levels are obtained by transferring (our chains consist solely of free Mackey functors). Secondly, the bottom level tensor product itself gets arbitrarily large as we increase $n,m$, making it impractical to compute with it. The extra complexity is reflected in the fact that the results for representations $m\lambda-n\sigma$ and $n\sigma-m\lambda$ (subsections \ref{Lambda minus sigma even} to \ref{Sigma minus Lambda odd}) are more involved than those for $n\sigma+m\lambda$ and $-n\sigma-m\lambda$ (subsections \ref{PureStart} to \ref{PureEnd}); after all, we need not use any box products for the representations of the form $n\sigma+m\lambda$ and $-n\sigma-m\lambda$.

While computing the box product of these chain complexes by hand is very impractical, a computer can do it efficiently. The idea is that our chains consist of solely free Mackey functors over $\underline{\Z}$, so every differential can be completely described by a matrix with integer entries. The operations of transfer, restriction and group action can all be performed algorithmically for free Mackey functors, and their effect can be described in terms of these matrices. Similarly, the tensor product can also be computed algorithmically, and then the box product is just obtained by transferring it to higher levels. At the final step, we need to take homology and that can be achieved via a Smith Normal Form algorithm over $\Z$.

There are a few more technicalities in this procedure that we haven't addressed here, but once these details are dealt with, this process allows us to algorithmically compute the additive structure of the $RO(G)$ homology, in any given range for our representations (for $G=C_4$, this amounts to a given range for $n,m$). 

For the multiplicative structure we need to be able to compute the product of any two generators. Just like with tensor products, this can be directly performed only on the bottom level. If the generators live in a higher level, the idea is to first restrict them to the bottom level, multiply these restrictions, and then invert the restriction map. This is possible because in free Mackey functors, restrictions are injective, and our chain complexes consist exclusively of such Mackey functors. 

So far we have enough information to \emph{verify} the multiplicative structure as it appears in section \ref{Results}, but not automatically \emph{compute} it. In other words, the expressions of the generators in section \ref{Results} need to be known a-priori and then the program can prove them in a user-specified range. But there is a final algorithm that eliminates this need, and allows us to automatically write our generators in terms of Euler and orientation classes. This "factorization" algorithm works by forming a multiplication table for the $RO(G)$ homology, and then turns this table into a colored graph, somewhat analogous to the Cayley graph of a group. There are two colors, corresponding to multiplication and division, and traversing this graph is equivalent to generating expressions of the generators like those appearing in section \ref{Results}.

This chains-based approach also works remarkably well with Massey products. And indeed, our program can compute Massey products, and their indeterminacy, directly from their definition. Finally, we can replace $\underline{\mathbb Z}$ with other constant Green functors such as $\underline{\mathbb F_p}$ for prime $p$ and $\underline{\mathbb Q}$.

In version 3.0 of the program, we introduced support for computing the $RO(G)$ graded homology of $G$-spaces other than the point. This extension is used in \cite{BC4S2} for the $C_4$ classifying space of the group $\Sigma_2$.

Since version 1.0, we have significantly improved the program's runtime performance and memory usage. This was partly achieved by using the sparse matrix format to store the differentials; in this context, a matrix is sparse if the vast majority of its entries are $0$ and thus we need not waste memory storing them. Another improvement was made by introducing a variant of algebraic Morse theory (see \cite{Lam19}) that preserves equivariance. This allows us to reduce our chain complexes to smaller ones in the same equivariant homotopy type. With this reduction, we can compress box products of chain complexes anywhere between 30\% to 90\%, with larger chain complexes leading to better compression ratios.

All these details plus many more can be found in the documentation for our code available \href{https://nickg-math.github.io/Mackey/html/index.html}{here}. 

The reader only interested in testing our program (in the case $G=C_4$) can simply download the executable for their operating system available \href{https://github.com/NickG-Math/Mackey/releases}{here}; no programming knowledge is required to run it.

The source code itself is written in C++ and hosted on a \href{https://github.com/NickG-Math/Mackey.git}{Github repository} to encourage participation and contribution. We have tried to make the code modular and extensible while at the same time fully \href{https://nickg-math.github.io/Mackey/html/index.html}{documenting} both how to use it and how it works under the hood.

\section{The Mackey functors}\label{MackeyFunctors}

The data in a $C_4$-Mackey functor $\underline M$ can be depicted using a Lewis diagram
\begin{equation}
\begin{tikzcd}
\underline{M}(C_4/C_4)\ar[d, "\Res^4_2" left, bend right]\\
\underline{M}(C_4/C_2)\ar[u, "\Tr_2^4" right,bend right]\ar[d, "\Res^2_1" left, bend right]\ar[loop right,distance=3em, "C_4/C_2"]\\
\underline{M}(C_4/e)\ar[u, "\Tr_1^2" right,bend right]\ar[loop right,distance=3em, "C_4"]
\end{tikzcd}
\end{equation}
We shall refer to $\underline{M}(C_4/C_4)$ as the top level (or $C_4$ level), to $\underline{M}(C_4/C_2)$ as the middle level (or $C_2$ level) and finally to $\underline{M}(C_4/e)$ as the bottom level (or $e$ level).\smallbreak

To improve readability, we shall stop underlining our Mackey functors. The only potential point of confusion is $\Z$ which could either denote the trivial $C_4$ module or the fixed point Mackey functor associated to it. Which one we mean will usually be clear from the context, but when the distinction is important we shall underline the Mackey functor $\underline{\Z}$.

We will also write $H_*(-)$ in place of $\underline{H}_*^{C_4}(-;\underline{\Z})$ and $H_{\bigstar}(-)$ for $\underline{H}_{\bigstar}^{C_4}(-;\underline{\Z})$; the little and big asterisks stand for integer and $RO(C_4)$ grading respectively. \smallbreak

The real representation ring $RO(C_4)$ is generated by the irreducible representations $\sigma$ and $\lambda$ where $\sigma=\R$ is reflection and $\lambda=\R^2$ is rotation by $\pi/2$, both leaving $0$ fixed. So the computation of $H_{\bigstar}(S)$ breaks down to calculating $H_*S^{\pm n\sigma\pm m\lambda}$ for the four possible sign combinations. Here and throughout this paper, $n,m$ will always stand for nonnegative integers.
\medbreak

We now display the Lewis diagrams of the Mackey functors appearing in our computations. 

For $H_*S^{n\sigma+m\lambda}$ we have the 5 Mackey functors:
\begin{equation}
\underline{\Z}=\begin{tikzcd}
\Z\ar[d, "1" left, bend right]\\
\Z\ar[u, "2" right,bend right]\ar[d, "1" left, bend right]\\
\Z\ar[u, "2" right,bend right]
\end{tikzcd}\quad \quad \quad
\Z_{-}=\begin{tikzcd}
0\ar[d, bend right]\\
\Z\ar[u,bend right]\ar[d, "1" left, bend right]\ar[loop right,"-1"]\\
\Z\ar[u, "2" right,bend right]\ar[loop right,"-1"]
\end{tikzcd}
\end{equation}
	\begin{equation}
\ev{\Z/4}=\begin{tikzcd}
\Z/4\ar[d, "1" left, bend right]\\
\Z/2\ar[u, "2" right,bend right]\ar[d, bend right]\\
0\ar[u, bend right]
\end{tikzcd}\quad \quad \quad
\ev{\Z/2}=\begin{tikzcd}
\Z/2\ar[d, bend right]\\
0\ar[u,bend right]\ar[d, bend right]\\
0\ar[u, bend right]
\end{tikzcd}
\quad \quad \quad
\overline{\ev{\Z/2}}=\begin{tikzcd}
0\ar[d, bend right]\\
\Z/2\ar[u,bend right]\ar[d, bend right]\\
0\ar[u, bend right]
\end{tikzcd}
\end{equation}

For $H_*S^{-n\sigma-m\lambda}$ we have the 4 additional Mackey functors:
\begin{gather}
L=\begin{tikzcd}
\Z\ar[d, "2" left, bend right]\\
\Z\ar[u, "1" right,bend right]\ar[d, "2" left, bend right]\\
\Z\ar[u, "1" right,bend right]
\end{tikzcd}\quad \quad \quad
p^*L=\begin{tikzcd}
\Z\ar[d, "2" left, bend right]\\
\Z\ar[u, "1" right,bend right]\ar[d, "1" left, bend right]\\
\Z\ar[u, "2" right,bend right]
\end{tikzcd}\quad \quad \quad
L_{-}=\begin{tikzcd}
\Z/2\ar[d, left, bend right]\\
\Z\ar[u, "1" right,bend right]\ar[d, "2" left, bend right]\ar[loop right,"-1"]\\
\Z\ar[u, "1" right,bend right]\ar[loop right,"-1"]
\end{tikzcd}\quad \quad \quad
p^*L_{-}=\begin{tikzcd}
\Z/2\ar[d, left, bend right]\\
\Z\ar[u, "1" right,bend right]\ar[d, "1" left, bend right]\ar[loop right,"-1"]\\
\Z\ar[u, "2" right,bend right]\ar[loop right,"-1"]
\end{tikzcd}	\end{gather}
Here $p^*$ denotes the functor from $C_2$ Mackey functors to $C_4$ Mackey functors induced by the quotient map $p:C_4\to C_4/C_2$.

For $H_*S^{m\lambda-n\sigma}$ we also have the trivial extension $\ev{\Z/2}\oplus \overline{\ev{\Z/2}}$ and the Mackey functor
\begin{equation}
Q=\begin{tikzcd}
\Z/2\ar[d, "0" left, bend right]\\
\Z/2\ar[u, "1" right,bend right]\ar[d, bend right]\\
0\ar[u,bend right]
\end{tikzcd}	\end{equation}
We have the nontrivial extensions
\begin{gather}
0\to \ev{\Z/2}\to Q\to \overline{\ev{\Z/2}}\to 0\\
0\to Q\to \ev{\Z/4}\to \ev{\Z/2}\to 0\end{gather}	

\medbreak

The additional Mackey functors present in $H_*S^{n\sigma-m\lambda}$ are the trivial extensions $\Z_{-}\oplus \ev{\Z/2}$,  $L\oplus \ev{\Z/2}$ and the 3 Mackey functors:
	\begin{equation}
L^{\sharp}=\begin{tikzcd}
\Z\ar[d, "1" left, bend right]\\
\Z\ar[u, "2" right,bend right]\ar[d, "2" left, bend right]\\
\Z\ar[u, "1" right, bend right]
\end{tikzcd}\quad \quad \quad
Q^{\sharp}=\begin{tikzcd}
\Z/2\ar[d, "1" left, bend right]\\
\Z/2\ar[u, "0" right, bend right]\ar[d, left, bend right]\\
0\ar[u, right, bend right]
\end{tikzcd}	\quad \quad \quad
	\Z_{-}^{\flat}=\begin{tikzcd}
0\ar[d, bend right]\\
\Z\ar[u, bend right]\ar[d, "2" left, bend right]\ar[loop right, "-1"]\\
\Z\ar[u, "1" right, bend right]\ar[loop right, "-1"]
\end{tikzcd}
\end{equation} 
The sharp operation $\sharp$ exchanges $\Res^4_2$ and $\Tr^4_2$ in our Mackey functor, while the flat operation $\flat$ exchanges $\Res^2_1$ and $\Tr^2_1$. For example $p^*L=L^{\flat}$ and $p^*L_{-}=L_{-}^{\flat}$.\medbreak

We have the nontrivial extensions
\begin{gather}
0\to  \ev{\Z/2}\to  \ev{\Z/4}\to Q^{\sharp}\to 0\\
0\to  \overline{\ev{\Z/2}}\to  Q^{\sharp}\to \ev{\Z/2}\to 0\\
0\to L\to L^{\sharp}\to \ev{\Z/2}\to 0\\
0\to \ev{\Z/2}\to L_{-}\to \Z_{-}^{\flat}\to 0
\end{gather}

\section{The generators}\label{Generators}	

As a $C_4$-Mackey functor $M$ has three levels, a generator for $M$ consists of three elements $a,b,c$ that generate the abelian groups $M(C_4/C_4), M(C_4/C_2)$ and $M(C_4/e)$ respectively. We shall employ the notation 
\begin{equation}
a|b|c
\end{equation}
to denote the top, middle and bottom generators in this order. 

Now $\pi_{\bigstar}(H\underline{\Z})$ is not just a (graded) Mackey functor, as it has a multiplicative structure making it a (graded) Green functor. Multiplication is performed levelwise and the Frobenius relation holds:
\begin{equation}
\Tr_K^H(x\Res^H_Ky)=\Tr_K^H(x)y
\end{equation}
where $K\subseteq H$ are subgroups of $C_4$.

 In this section we shall expound on the interplay between the Mackey functor and multiplicative structures, and demonstrate how every generator can be written in terms of the Euler and orientation classes. We begin by defining these classes in greater generality, following \cite{HHR16}.\medbreak

For any real representation $V$ of a group $G$ we have the Euler class
\begin{equation}
a_V:S^0\to S^V
\end{equation}
given by the inclusion of the north and south poles $0,\infty$. We shall only consider the image of $a_V$ in homology but it's important for some arguments to note that $a_V$ is defined on the sphere level.

If $V$ is orientable, namely the map $G\to GL(n)$ defining $V$ has positive determinant, then we have
\begin{equation}
H_{n}^G(S^V;\underline{\Z})=\underline{\Z}
\end{equation} 
(cf \cite{HHR16}). Choosing an orientation for $V$ determines an orientation class $u_V$ as the generator of the top level of this $\underline{\Z}$. Without orienting $V$ there is a sign ambiguity for $u_V$.

In \cite{HHR17} the following properties are proven in the case of $G=C_{2^n}$ (whenever $u_V$ appears it is implicit that $V$ is oriented).\medbreak

\begin{itemize}
	\item $a_Va_W=a_{V+W}$ and $u_{V+W}=u_Vu_W$
	\item $\Res^G_Ha_V=a_{\Res^G_HV}$ and $\Res^G_Hu_V=u_{\Res^G_HV}$ 
%	\item $N^G_Ha_V=a_{N^G_HV}$ and $u_{N_H^G|V|}N_H^G(u_V)=u_{N_H^GV}$ where norm of representations is induction, and $|V|$ denotes the trivial representation of dimension that of $V$.
	\item $|G:Stab(V)|a_V=0$ where $Stab(V)$ is the stabilizer (isotropy subgroup) of $V$.
	\item The Gold (au) Relation: If $V,W$ have dimension $2$ and $Stab(V)\le Stab(W)$,
	\begin{equation}
	a_Wu_V=	|Stab(W):Stab(V)|\cdot a_Vu_W
	\end{equation}
\end{itemize}

In our case, the real $C_4$ representations are spanned by $1,\sigma,\lambda$ and the orientable ones are spanned by $1,2\sigma,\lambda$. Therefore we have the classes
\begin{equation}
a_{\sigma},a_{\lambda},u_{2\sigma},u_{\lambda}
\end{equation}
living in the top level of $H_{\bigstar}S$. While $\sigma$ is not orientable as a $C_4$ representation, its restriction to $C_2\subseteq C_4$ is the trivial $C_2$ representation so we can consider
\begin{equation}
u_{\sigma}
\end{equation}
living in the middle level. We choose orientations coherently so that
\begin{equation}
\Res^4_2(u_{2\sigma})=u_{\sigma}^2
\end{equation}
\medbreak

To simplify the notation, for an element $a$ living in some level of a Mackey functor, we shall write $\bar a$ for its restriction to the level directly below. If $a$ lives in the top level, we can restrict $\bar a$ again and then $\bar{\bar{a}}$ will be the restriction two levels down. This notation is consistent with \cite{HHR17}. \medbreak

The Euler and orientation classes generate the following Mackey functors:
\begin{align}
a_{\sigma}|\zero|\zero&\rightsquigarrow \ev{\Z/2}\\
a_{\lambda}|\bar a_{\lambda}|\zero&\rightsquigarrow \ev{\Z/4}\\
u_{2\sigma}|u_{\sigma}^2|\bar u_{\sigma}^2&\rightsquigarrow \Z\\
u_{\lambda}|\bar u_{\lambda}|\bar{\bar{u}}_{\lambda}&\rightsquigarrow \Z\\
\zero|u_{\sigma}|\bar u_{\sigma}&\rightsquigarrow \Z_{-}
\end{align}
The Mackey functors themselves imply relations on these classes eg $2a_{\sigma}=0$. Moreover, since $H_{\bigstar}(S)$ is a Green functor we also have the Frobenius relation:
\begin{equation}
\Tr_K^H(x\Res^H_Ky)=\Tr_K^H(x)y
\end{equation}
We will refer to all these as secondary relations; the primary relations are those not implied by the additive (Mackey functor) structure or Frobenius. The only primary relation we have so far is the Gold relation:
\begin{equation}
a_{\sigma}^2u_{\lambda}=2u_{2\sigma}a_{\lambda}
\end{equation}
The Euler and orientation classes generate multiplicatively all of $H_*(S^{n\sigma+m\lambda})$.\medbreak

Before we explain how $H_*(S^{-n\sigma-m\lambda})$ is generated we need to take a moment and clarify what we mean by division:

\subsection{A digression on divisibilities}\label{Division}

Suppose we have elements $x\in H_VS$ and $y\in H_WS$ living on the same level and that are not both $0$. We will say that $y/x$ exists if $H_{W-V}S$ has a cyclic subgroup $C$ such that multiplication by $x$ maps $C\subseteq H_{W-V}S$ isomorphically onto the cyclic subgroup $\ev{y}\subseteq H_WS$ generated by $y$:
	\begin{center}\begin{tikzcd}
		H_{W-V}S\ar[r,"x"]& H_WS\\
		C\ar[u,hook]\ar[r,"x" above, "\approx" below]&\ev{y}\ar[u,hook]
		\end{tikzcd}\end{center}
	If $C$ is unique with this property, then the preimage of $y$ under multiplication by $x$ is a single element in $H_{W-V}S$ denoted by $y/x$.

For example, $1/x$ exists iff $x$ is invertible, and in that case $1/x=x^{-1}$ (and we will continue to use the $x^{-1}$ notation for inverses).

However in general, $y/x$ is less ambiguous than $y^{-1}x$ as the latter notation might suggest that $y^{-1}$ exists by itself and is multiplied with $x$. For instance, $2/u_{2\sigma}$ exists because $H_{-2}(S^{-2\sigma})\xrightarrow{u_{2\sigma}}H_0(S)=\Z$ is an isomorphism onto $2\Z\subseteq \Z$ in the top level. On the other hand, $1/u_{2\sigma}$ does not exist.

Let us note here that if the subgroup $C$ in the definition above is not unique, then there are multiple candidates for $y/x$. We explain what to do in these cases, and exactly how it comes up in the $RO(C_4)$ homology of a point, in subsection \ref{Subtle}.\medbreak

Getting back to  $H_*(S^{-n\sigma-m\lambda})$, we will prove that $u_{\sigma}^{-1}$ and $\bar{\bar{u}}_{\lambda}^{-1}$ both exist. In fact, the following elements all exist:
\begin{equation}2/u_{2\sigma}^n, \quad 2/\bar u_{\lambda}^m, \quad  4/u_{\lambda}^m\text{ , } \quad  4/(u_{2\sigma}^nu_{\lambda}^m)\end{equation}
Now for odd $n\ge 3$ set
\begin{equation}
w_n=\Tr_2^4(u_{\sigma}^{-n})
\end{equation}
We don't consider $w_1$ because $\Tr_2^4(u_{\sigma}^{-1})=0$. Next for odd $n\ge 1$ and $m\ge 1$ set
\begin{equation}
x_{n,m}=\Tr_1^4(\bar u_{\sigma}^{-n}\bar{\bar{u}}_{\lambda}^{-m})
\end{equation}
The $w_n,x_{n,m}$ are all $2$-torsion elements and we have the divisibilities:
\begin{equation}
w_n/(a_{\sigma}^{i}a_{\lambda}^{j}u_{\lambda}^{k}), \quad  x_{n,m}/a_{\sigma}^i
\end{equation}
\medbreak

The first element not obtained by Euler and orientation classes through the operations of multiplication, division (wherever possible), transfers and restrictions is the generator $s$ in the top level of $H_{-3}S^{-2\lambda}=\ev{\Z/4}$. We have the divisibilities
\begin{equation}
s/(u_{2\sigma}^{i}a_{\lambda}^{j}u_{\lambda}^{k}), \quad \bar s/(u_{\sigma}^i\bar a_{\lambda}^{j}\bar u_{\lambda}^{k})
\end{equation}
Thus far we have accounted for every element in $H_*(S^{n\sigma+m\lambda})$ and $H_*(S^{-n\sigma-m\lambda})$.\medbreak

For $H_*(S^{m\lambda-n\sigma})$ we have additional elements
\begin{equation}
u_{\lambda}/u_{2\sigma}^{i}, \quad  (2a_{\lambda})/(a_{\sigma}u_{2\sigma}^i)\end{equation}
and for  $H_*(S^{n\sigma-m\lambda})$ we have
\begin{equation}
(2u_{2\sigma})/u_{\lambda}, \quad (4u_{2\sigma})/u_{\lambda}^i, \quad a_{\sigma}^2/a_{\lambda},\quad a_{\sigma}^3/a_{\lambda}^{m}
\end{equation}
We also obtain the relations:
\begin{gather}
2s=w_3(a_{\sigma}^3/a_{\lambda}^2)\\
a_{\sigma}s=\Tr_2^4((2u_{\sigma})/\bar u_{\lambda}^2)
\end{gather}
In the second equation, multiplication by $a_{\sigma}$ is the projection $\Z/4\to \Z/2$ so we equivalently have
\begin{gather}
2s=w_3(a_{\sigma}^3/a_{\lambda}^2)\\
s\text{ mod }2=\Tr_2^4((2u_{\sigma})/\bar u_{\lambda}^2)/a_{\sigma}
\end{gather}
expressing $2s$ and $s$ mod $2$ in terms of Euler and orientation classes. Thus $s$ is obtained from Euler and orientation classes through the extension
\begin{equation}
0\to \Z/2\to \Z/4\to \Z/2\to 0
\end{equation}
(note: the extension determines $s$ only up to a sign; i.e. $s$ cannot be canonically chosen from this extension). But if we want $\Res^4_2(s)$, then we need to replace this group extension with one of Mackey functors. In fact, in this case we have two such extensions:
\begin{gather}
0\to \ev{\Z/2}\to L_{-}\to \Z_{-}^{\flat}\to 0\\
0\to \overline{\ev{\Z/2}}\to Q^{\sharp}\to \ev{\Z/2}\to 0
\end{gather}
\begin{itemize}
	\item In the first extension, $a_{\sigma}s|\zero|\zero$ generates $\ev{\Z/2}$ and  $\zero|(2u_{\sigma})/\bar u_{\lambda}^2|\bar u_{\sigma}\bar {\bar u}_{\lambda}^{-2}$ generates $\Z_{-}^{\flat}$.
	\item In the second extension, $\zero|u_{\sigma}^3\bar s|\zero$ generates $\overline{\ev{\Z/2}}$ and $a_{\sigma}^3/a_{\lambda}^2|\zero|\zero$ generates $\ev{\Z/2}$.
\end{itemize} 

From this description of the generators, we see that first short exact sequence gives the formula for $a_{\sigma}s$ (which is equivalent to the formula for $s$ mod $2$), while the second gives
\begin{equation}
\bar s=\Res^4_2(a_{\sigma}^3/a_{\lambda}^{2})/u_{\sigma}^3
\end{equation}
Applying $\Tr_2^4$ on both sides returns the formula for $2s$.
\medbreak

In order to summarize this whole discussion more concisely, we will use a more general notion of "generator". In this notion, the span of a list of elements is not just polynomials on those generators combined with transfers and restrictions (the Green functor span) but we will also allow any divisibilities that occur as well as Mackey functor extensions in which the outer two Mackey functors are already in the span. For this to be well defined we need to note which divisibilities and extensions actually occur.

This generalized notion satisfies the following property: If we have two Green functor maps $f,g:M\to N$ and $M$ has a set $A$ of generalized generators then $f=g$ on $M$ iff $f=g$ on $A$ and $f=g$ on a generator for any extension that occurs (after all, these generators cannot be canonically chosen through the extensions).

The other part of this property has to do with whether or not a map $f:A\to N$ extends to a Green functor map $f:M\to N$. This is of course tantamount to $f$ satisfying all Green functor relations. Ideally we would like to only list the primary relations on the generalized generators $A$ and recover all other relations from these and the secondary ones (which result from the additive structure and the Frobenius relations). As we explain in subsection \ref{Subtle}, this might not always be possible. In the special case of $\underline \pi_{\bigstar}^{C_4}(H\underline{\Z})$ it does however work out so we can legitimately call them "primary relations".

With this language, we can summarize this section as follows: The Green functor $\underline \pi_{\bigstar}^{C_4}(H\underline{\Z})$ has generalized generators $a_{\sigma},u_{2\sigma},u_{\sigma},a_{\lambda},u_{\lambda}$ in degrees $\bigstar=-\sigma,2-2\sigma,1-\sigma,-\lambda,2-2\lambda$ respectively. These classes individually generate the Mackey functors 
\begin{align}
a_{\sigma}|\zero|\zero&\rightsquigarrow \ev{\Z/2}\\
a_{\lambda}|\bar a_{\lambda}|\zero&\rightsquigarrow \ev{\Z/4}\\
u_{2\sigma}|u_{\sigma}^2|\bar u_{\sigma}^2&\rightsquigarrow \Z\\
u_{\lambda}|\bar u_{\lambda}|\bar{\bar{u}}_{\lambda}&\rightsquigarrow \Z\\
\zero|u_{\sigma}|\bar u_{\sigma}&\rightsquigarrow \Z_{-}
\end{align}
 and the only primary relation is the Gold relation
\begin{gather}
a_{\sigma}^2u_{\lambda}=2u_{2\sigma}a_{\lambda}\end{gather}
The only extension that occurs is for the generator $s|\bar s|\zero$ of the $\ev{\Z/4}$ in dimension $\bigstar=-3+2\lambda$ and is specified by:
\begin{gather}
\bar s=\Res^4_2(a_{\sigma}^3/a_{\lambda}^{2})/u_{\sigma}^3\\
a_{\sigma}s=\Tr_2^4((2u_{\sigma})/\bar u_{\lambda}^2)
\end{gather}

There are many divisibilities that occur and we have indicated most of them earlier in this section; an exhaustive list is included in the next section together with the complete determination of the additive structure.\medbreak

\subsection{Some technical remarks}\label{Subtle} We end this section with a few subtle points that can arise when dealing with quotients $y/x$.

First, there might be multiple choices of $y/x$ when it exists (cf subsection \ref{Division}). That is, there can be multiple cyclic subgroups $C$ with $C\xrightarrow{x}\ev{y}$ an isomorphism. In that case, we should choose $y/x$ so that
\begin{equation}
z\cdot (y/(xz))=y/x
\end{equation}
for any $z$ such that $y/(xz)$ exists. If there are multiple candidates for $y/(xz)$ for some $z$, we have to ensure coherency in our choices so that not only the above equation is true, but also
\begin{equation}
w\cdot (y/(xzw))=y/(xz)
\end{equation}
for any $w$ such that $y/(xzw)$ exists and so on.

In practice, the different candidates can be distinguished by their products with Euler and orientation classes. For example, there may only be one candidate whose product with $a_{\sigma}$ is zero, i.e. a transfer; if we make that choice then coherence is guaranteed by the Frobenius relation.\medbreak

In the case of the $RO(C_4)$ homology of a point, the only cases where this situation can arise have to do with the elements
\begin{equation}
\Tr_2^4\left((2u_{\sigma}^{n})/\bar u_{\lambda}^m\right)\text{ , }(a_{\sigma}^{2m}u_{2\sigma}^{n/2-m})/a_{\lambda}^m
\end{equation}
spanning $H_{n-2m}=\Z\oplus \Z/2$ for $n-2m\ge 0$ and $m\ge 2$. Multiplication by $u_{\lambda}^m$ maps the former to $4u_{2\sigma}^{n/2}$ and the latter to $0$, hence the two choices for $(4u_{2\sigma}^{n/2})/u_{\lambda}^m$ are:
\begin{equation}
\Tr_2^4\left((2u_{\sigma}^{n})/\bar u_{\lambda}^m\right)\text{ , }\Tr_2^4\left((2u_{\sigma}^{n})/\bar u_{\lambda}^m\right)+(a_{\sigma}^{2m}u_{2\sigma}^{n/2-m})/a_{\lambda}^m
\end{equation}
The second element is not a transfer, so we pick the first:
\begin{equation}
(4u_{2\sigma}^{n/2})/u_{\lambda}^m=\Tr_2^4\left((2u_{\sigma}^{n})/\bar u_{\lambda}^m\right)
\end{equation}
There is another benefit to this choice: For $n=2$ there is a unique candidate for $(4u_{2\sigma})/u_{\lambda}^m$, and our general choice satisfies the nice property:
\begin{equation}
(4u_{2\sigma}^{n/2})/u_{\lambda}^m=((4u_{2\sigma})/u_{\lambda}^m)\cdot u_{2\sigma}^{n/2-1}
\end{equation}

This brings us to the second subtle point: The expressions $(x/z)\cdot (y/w)$ and $(xy)/(zw)$ are not always equivalent: one can exist when the other doesn't, and even if both exist then they might not be equal! Case in point:
\begin{equation}
w_3(a_{\sigma}^3/a_{\lambda}^2)\neq (w_3a_{\sigma}^3)/a_{\lambda}^2
\end{equation}
as the left hand side generates a $\Z/2$, while the right is trivial owing to $w_3a_{\sigma}^3=0$.

If $x/z,y/w$ exist then $(x/z)\cdot (y/w)=(xy)/(zw)$ is equivalent to $(x/z)\cdot (y/w)$ and $xy$ generating isomorphic cyclic subgroups. In practice, the elements given by our spectral sequences are of the form $(x/z)\cdot (y/w)$ and the additive structure is known apriori, so we can readily check this equality.

For the generators displayed in section \ref{Results}, the only instance where  $(x/z)\cdot (y/w)$ and $(xy)/(zw)$ differ happens with
\begin{equation}
u_{2\sigma}^i(s/(a_{\lambda}^ju_{\lambda}^k))
\end{equation}
for $i,j,k\ge 0$ and $i,k>0$. This element generates a $\Z/4$ while $u_{2\sigma}^is$ generates a $\Z/2$ and thus $(u_{2\sigma}^is)/(a_{\lambda}^ju_{\lambda}^k)$ does \emph{not} equal the $\Z/4$ generator, but is rather the mod $2$ reduction of that generator.

There's another problem that stems from this point and it has to do with relations. We want to be able to reduce relations on $x/y$ to equivalent relations on $x$. A relation on $x/y$ takes the form $(x/y)\cdot z=0\in H_VS$ for some element $z\in H_{\bigstar}S$. If multiplication by $y$ is an isomorphism in $H_VS$ then we can clear denominators with $y$ and get the equivalent relation $xz=0$. If it's not an isomorphism then $xz=0$ may not be equivalent to $(x/y)\cdot z=0$.

Here's an example arising in "nature": Let's take the generator 
\begin{equation}y=(w_3a_{\lambda})/(a_{\sigma}u_{2\sigma})
\end{equation}
of a $\Z/2$. We want to establish the relation $u_{\lambda}y=0$. First, the homology group $u_{\lambda}y$ lives in is a $\Z/2$ (not $0$), and neither $u_{\lambda}$ nor $y$ are transfers (so we can't use the Frobenius relation); this means that $u_{\lambda}y=0$ is not a secondary relation. Second, $u_{\lambda}y=0$ is not equivalent to  $(u_{\lambda}y)a_{\sigma}u_{2\sigma}=0$ because multiplication by $a_{\sigma}u_{2\sigma}$ is not an isomorphism for the homology group $u_{\lambda}y$ lives in. 

However, the homology group $u_{\lambda}y$ lives in is generated by $(w_3a_{\lambda}^2)/a_{\sigma}^3$, as we can see from the tables in section \ref{Results}. So while multiplication by $a_{\sigma}u_{2\sigma}$ is not an isomorphism there, multiplication by $a_{\sigma}^3$ \emph{is}. Thus $u_{\lambda}y=0$ is equivalent to $a_{\sigma}^3u_{\lambda}y=0$ which is true because $a_{\sigma}^3u_{\lambda}=0$ by the Gold relation.

We employ a similar strategy when the homology group in the degree of the product is not cyclic. Let's take for example the relation
\begin{equation}
\left((2u_{2\sigma})/u_{\lambda}\right)^2=((4u_{2\sigma}^2)/u_{\lambda}^2)+a_{\sigma}^4/a_{\lambda}^2
\end{equation}
This "exotic multiplication" was pointed out to us by Mingcong Zeng. To prove it, write
\begin{equation}
\left((2u_{2\sigma})/u_{\lambda}\right)^2=x\cdot (4u_{2\sigma}^2)/u_{\lambda}^2+y\cdot a_{\sigma}^4/a_{\lambda}^2
\end{equation}
for unknown integers $x,y$. To find $x$ we can multiply by $u_{\lambda}^2$ and use the relation 
\begin{equation}
u_{\lambda}\cdot (a_{\sigma}^4/a_{\lambda}^2)=0
\end{equation}
proven in the Appendix \ref{Appendix}. To find $y$ we instead multiply with $a_{\lambda}^2$ and use 
\begin{equation}
a_{\lambda}\cdot (4u_{2\sigma}^2/u_{\lambda}^2)=0
\end{equation}
which is proven by appealing to Frobenius ($ (4u_{2\sigma}^2/u_{\lambda}^2)=\Tr_2^4((2u_{\sigma}^4)/\bar u_{\lambda}^2)$). In the end we get $x=y=1$ as desired.

This strategy fails when we have generators that are not of the form $x/y$. For example the seemingly innocuous relation
\begin{equation}
s\cdot ((2u_{2\sigma})/u_{\lambda})=2u_{2\sigma}(s/u_{\lambda})
\end{equation}
can't be proven by multiplying with $u_{\lambda}$ as the multiplication map can't distinguish between $0$ and $2u_{2\sigma}(s/u_{\lambda})$ due to $2su_{2\sigma}=0$. It's also not the image of another relation under multiplication by $u_{2\sigma}$. Instead, we can immediately deduce the relation from the general simple fact of denominator exchange:
\begin{equation}
(x/z)\cdot (y/w)=(x/w)\cdot (y/z)
\end{equation}
as long as $x/(zw),y/(zw)$ exist.\medbreak

It turns out that for the integral $RO(C_4)$ homology of a point, all relations can be recovered from the Gold and the secondary relations using the ideas above. Proving this is quite tedious, as we need to consider all unordered pairs of generators that are not transfers and compute their product in terms of the other generators. This work is displayed in the Appendix \ref{Appendix}.

\section{{\texorpdfstring{A conjecture for $G=C_{2^n}$}{A conjecture for G=C2n}}}\label{Conjecture}

For $G=C_{2^n}$ let $\sigma$ denote the sign representation and $\lambda_{k}$ denote the $2$-dimen\-sional representation given by rotation by $\pi/2^{k-1}$ degrees for $k=2,...,n$. Then the Mackey functor $\underline H_{-3}S^{-2\lambda_n}$ is on each orbit
\begin{equation}
\underline H_{-3}S^{-2\lambda_n}(G/C_{2^k})=\Z/2^k
\end{equation}
for $k>0$ and $\underline H_{-3}S^{-2\lambda_n}(G/e)=0$. Transfers are the usual inclusion maps $\Z/2^k\hookrightarrow \Z/2^l$ for $k\le l$, while restrictions are the projection maps $\Z/2^l\to \Z/2^k$ for $k\le l$.

We let $s_n$ denote a generator of $\Z/2^n$. Then we can directly compute that
\begin{equation}
s_na_{\sigma}=\Tr_{1}^{2^n}[(\Res^{2^{n-1}}_1u_{\sigma})(\Res^{2^n}_1 u_{\lambda_n})^{-2}]
\end{equation}
generating a $\Z/2$.

Now note that $2s_n$ is the transfer of the $C_{2^{n-1}}$ generator $s_{n-1}$ hence by induction, $s_n$ is generated by the Euler and orientation classes of $C_{2^k}$ for $2\le k\le n$ through the extension
\begin{equation}
0\to \Z/2^{n-1}(2s_n)\to \Z/2^n(s_n)\to \Z/2(s_na_{\sigma})\to 0
\end{equation}

\begin{conj}For all $G=C_{2^n}$, the Euler classes, orientation classes and $s_n$ together generate $\underline H_{\bigstar}S$ under the operations of multiplication, division, transfer and restriction.
\end{conj}
This has been verified in a finite range for $n\le 5$.

We further expect the Gold relation to generate all relations as in the preceding section. Finally, we make no conjecture for the additive structure, as computational data suggest that it's very complicated with no visible patterns.
\newpage
\section{The results}\label{Results}

We compile the results of our computation of the Green functor $H_*S^{\pm n\sigma\pm m\lambda}=\underline H_*^{C_4}(S^{\pm n\sigma\pm m\lambda};\underline{\Z})$ where $n,m\ge 0$ as always. We consider 8 separate cases based on the signs in $\pm n\sigma\pm m\lambda$ and the parity of $n$ (this parity determines whether the representation $\pm n\sigma\pm m\lambda$ is orientable or not); each case gets its own subsection containing both the additive and multiplicative structures. The 8 cases are ordered roughly in increasing complexity, which also happens to be the order in which we prove all these results in sections \ref{ProofsPureHomology} through \ref{Proofslambda}.

 The notation for the Mackey functors and their generators has been explained in the preceding two sections. For improved formatting we shall write $\frac{x}{y}$ in place of $x/y$.

\phantom{}\bigbreak

\subsection{\texorpdfstring{$\boldsymbol{H_*S^{n\sigma+m\lambda}}$ for even $\boldsymbol n$}{Sigma plus Lambda oriented}}\label{PureStart}

\begin{equation}
H_*(S^{n\sigma+m\lambda})=
\begin{cases}\hspace{-0.3em}
\begin{tabular}{l p{0.5em} l p{1.5em} l}
$\Z$			& if &  $*=n+2m$	&&\\
$\ev{\Z/4}$		& if & $n\le *<n+2m$		&and &$*$ is even\\
$\ev{\Z/2}$& if &  $0\le *<n$& and& $*$ is even\\
\end{tabular}
\end{cases}
\end{equation}

\def\arraystretch{2}
\setlength{\tabcolsep}{4pt}
\begin{center}
\begin{tabular}{p{0.2em} l p{4em} l p{1.5em} l}

	$\bullet$& $u_{2\sigma}^{n/2}u_{\lambda}^m|u_{\sigma}^n\bar u_{\lambda}^m|\bar u_{\sigma}^n \bar{\bar{u}}_{\lambda}^m$		& generates &$H_{n+2m}=\Z$		&	 &\\

	$\bullet$&$u_{2\sigma}^{n/2}a_{\lambda}^{m-i}u_{\lambda}^{i}|u_{\sigma}^n\bar a_{\lambda}^{m-i}\bar u_{\lambda}^i|\zero$	& generates &$H_{n+2i}=\ev{\Z/4}$	&for & $0\le i<m$\\

	$\bullet$& $a_{\sigma}^{n-2i}u_{2\sigma}^ia_{\lambda}^m|\zero|\zero$ 	& generates &$H_{2i}=\ev{\Z/2}$ 	&for& $0\le i\le \frac n2-1$\\

\end{tabular}
\end{center}

\phantom{}\bigbreak
\phantom{}\medbreak

\subsection{\texorpdfstring{$\boldsymbol{H_*S^{n\sigma+m\lambda}}$ for odd $\boldsymbol{n}$}{Sigma plus Lambda nonoriented}}

\begin{equation}
H_*(S^{n\sigma+m\lambda})=
\begin{cases}\hspace{-0.3em}
\begin{tabular}{l p{0.5em} l p{1.5em} l}
$\Z_{-}$			& if &  $*=n+2m$	&&\\
$\overline{\ev{\Z/2}}$		& if &	$n\le *<n+2m$		&and &$*$ is odd\\
$\ev{\Z/2}$& if & $0\le *<n+2m$& and& $*$ is even\\  
\end{tabular}
\end{cases}
\end{equation}

\begin{center}
	\begin{tabular}{p{0.2em} l p{4em} l p{1.5em} l}
		
	$\bullet$&	$\zero|u_{\sigma}^n\bar u_{\lambda}^m|\bar u_{\sigma}^n\bar{\bar{u}}_{\lambda}^m$		& generates & $H_{n+2m}=\Z_{-}$		&	 &\\

	$\bullet$&	$\zero|u_{\sigma}^n\bar a_{\lambda}^{m-i}\bar u_{\lambda}^i|\zero$	& generates &$H_{n+2i}=\overline{\ev{\Z/2}}$	&for & $0\le i<m$\\
		
	$\bullet$&	$a_{\sigma}^{n-2i}u_{2\sigma}^ia_{\lambda}^m|\zero|\zero$	& generates &$H_{2i}=\ev{\Z/2}$ 	&for& $0\le i\le \frac n2-1$\\
		
	$\bullet$&	$a_{\sigma}u_{2\sigma}^{(n-1)/2}a_{\lambda}^{m-i}u_{\lambda}^i|\zero|\zero$& generates &$H_{n+2i-1}=\ev{\Z/2}$ 	&for& $1\le i<m$\\
	\end{tabular}
\end{center}

\newpage

\subsection{\texorpdfstring{$\boldsymbol{H_*S^{-n\sigma-m\lambda}}$ for even $\boldsymbol n$}{Minus Sigma minus Lambda oriented}}
If $n,m$ are not both $0$,

\begin{equation}
H_*(S^{-n\sigma-m\lambda})=
\begin{cases}
\hspace{-0.3em}
\begin{tabular}{l p{0.5em} l p{1.5em} l p{1.5em} l}
$L$			& if &  $*=-n-2m$		&and & $m\neq 0$&&\\
$p^*L$		& if &	$*=-n-2m$		&and & $m=0$&&\\
$\ev{\Z/4}$	& if &	$-n-2m<*<-n-1$	&and & $*$ is odd&&\\
$\ev{\Z/2}$& if & $-n-1\le *<-1$& and& $*$ is odd &and& $m\neq 0$\\
$\ev{\Z/2}$& if & $-n+1\le *<-1$& and& $*$ is odd  &and& $m=0$\\
\end{tabular}
\end{cases}
\end{equation}

\def\arraystretch{3}
\setlength{\tabcolsep}{4pt}
\begin{center}
	\begin{tabular}{p{0.2em} l p{4em} l p{1.5em} l}
	$\bullet$&	$\dfrac4{u_{2\sigma}^{n/2}u_{\lambda}^{m}}\Big\rvert \dfrac 2{u_{\sigma}^n\bar u_{\lambda}^m}\Big\rvert \dfrac1{\bar u_{\sigma}^{n}\bar{\bar{u}}_{\lambda}^{m}}$		&		 generates & $H_{-n-2m}=L$		 &for &$m\neq 0$\\

	$\bullet$&	$\dfrac{2}{u_{2\sigma}^{n/2}}\Big\rvert\dfrac{1}{u_{\sigma}^{n}}\Big\rvert\dfrac{1}{\bar u_{\sigma}^{n}}$  &generates	&$H_{-n}=p^*L$	 &for & $m=0$ \\
		
	$\bullet$&	$\dfrac{s}{u_{2\sigma}^{n/2}a_{\lambda}^{i-2}u_{\lambda}^{m-i}}\Big\rvert\dfrac{\bar s}{u_{\sigma}^{n}\bar a_{\lambda}^{i-2}\bar u_{\lambda}^{m-i}}\Big\rvert \zero$ & generates & $H_{-n-2m+2i-3}=\ev{\Z/4}$  &for & $2\le i\le m$\\
		
		$\bullet$&	$\dfrac{x_{2i+1,1}}{a_{\sigma}^{n-2i-1}a_{\lambda}^{m-1}}\Big\rvert \zero \Big\rvert \zero$ &generates &$H_{-2i-3}=\ev{\Z/2}$  &for & $0\le i\le\frac{n}2-1$, $m\neq 0$\\

$\bullet$&	$\dfrac{w_{2i+3}}{a_{\sigma}^{n-2i-3}}\Big\rvert \zero \Big\rvert \zero$ &generates &$H_{-2i-3}=\ev{\Z/2}$  &for & $0\le i<\frac{n}2-1$, $m=0$\\
	\end{tabular}
\end{center}

\newpage
\subsection{\texorpdfstring{$\boldsymbol{H_*S^{-n\sigma-m\lambda}}$ for odd $\boldsymbol n$}{Minus Sigma minus Lambda nonoriented}}\label{PureEnd}

\begin{equation}
H_*(S^{-n\sigma-m\lambda})=
\begin{cases}
\hspace{-0.3em}
\begin{tabular}{l p{0.5em} l p{1.5em} l}
$L_{-}$			& if &  $*=-n-2m$		&and & $m\neq 0$\\
$p^*L_{-}$		& if &	$*=-n$		&and & $n>1,m=0$\\
$\Z_{-}$		& if &	$*=-1$		&and & $n=1,m=0$\\
$\overline{\ev{\Z/2}}$	& if &	$-n-2m<*<-n-1$	&and & $*$ is even\\
$\ev{\Z/2}$& if & $-n-2m<*<-1$& and& $*$ is odd\\
\end{tabular}
\end{cases}
\end{equation}

\begin{center}
	\begin{tabular}{p{0.2em} l p{4em} l p{1.5em} l}
		$\bullet$&	$x_{n,m}\Big\rvert \dfrac 2{u_{\sigma}^n\bar u_{\lambda}^m}\Big\rvert \dfrac1{\bar u_{\sigma}^{n}\bar{\bar{u}}_{\lambda}^{m}}$		&		 generates & $H_{-n-2m}=L_{-}$		 &for &$m\neq 0$\\

		$\bullet$&	$w_n|u_{\sigma}^{-n}| \bar u_{\sigma}^{-n}$  &generates	&$H_{-n}=p^*L_{-}$	 &for & $n>1,m=0$ \\
		
		$\bullet$&	$\zero|u_{\sigma}^{-1}|\bar u_{\sigma}^{-1}$  &generates	&$H_{-1}=\Z_{-}$	 &for & $n=1$, $m=0$ \\

		$\bullet$&	$\zero\Big\rvert\dfrac{\bar s}{u_{\sigma}^{n}\bar a_{\lambda}^{i-2}\bar u_{\lambda}^{m-i}}\Big\rvert \zero$ & generates & $H_{-n-2m+2i-3}=\overline{\ev{\Z/2}}$  &for & $2\le i\le m$\\
		
		$\bullet$&	$\dfrac{2s}{a_{\sigma}u_{2\sigma}^{(n-1)/2}a_{\lambda}^{i-2}u_{\lambda}^{m-i}}\Big\rvert \zero\Big\rvert\zero$
		
		%$\dfrac{x_{n,i}}{a_{\lambda}^{m-i}}\Big\rvert \zero \Big\rvert \zero$
		
		 &generates &$H_{-n-2m+2i-2}=\ev{\Z/2}$  &for & $2\le  i\le m$\\
		
		$\bullet$&	$\dfrac{x_{2i+1,1}}{a_{\sigma}^{n-2i-1}a_{\lambda}^{m-1}}\Big\rvert \zero \Big\rvert \zero$ &generates &$H_{-2i-3}=\ev{\Z/2}$  &for & $0\le i<\frac{n-1}2$, $m\neq 0$\\
		
		$\bullet$&	$\dfrac{w_{2i+3}}{a_{\sigma}^{n-2i-3}}\Big\rvert \zero \Big\rvert \zero$ &generates &$H_{-2i-3}=\ev{\Z/2}$  &for & $0\le i<\frac{n-3}2$, $m=0$\\
	\end{tabular}
\end{center}

\newpage

\subsection{\texorpdfstring{$\boldsymbol{H_*S^{m\lambda-n\sigma}}$ for even $\boldsymbol n$}{Lambda minus Sigma oriented}}\label{Lambda minus sigma even} If $n,m$ are both nonzero,

\begin{equation}
H_*(S^{m\lambda-n\sigma})=
\begin{cases}
\hspace{-0.3em}
\begin{tabular}{l p{0.5em} l p{1.5em} l}
$\Z$			& if &  $*=2m-n$		& &\\
$\ev{\Z/4}$		& if &$-n+2\le *<2m-n$	&and &  $*$ is even\\
$\ev{\Z/2}$		& if &$-n+1\le *\le -3$ 		&and & $*$ is odd\\
$Q$				& if & $*=-n$&&
\end{tabular}
\end{cases}
\end{equation}
\begin{center}
	\begin{tabular}{p{0.2em} l p{4em} l p{1.5em} l}
		$\bullet$&	$\dfrac{u_{\lambda}^m}{u_{2\sigma}^{n/2}}\Big\rvert \dfrac{\bar u_{\lambda}^m}{u_{\sigma}^{n}}\Big\rvert \dfrac{\bar{\bar u}_{\lambda}^m}{\bar u_{\sigma}^n}$		&		 generates &  $H_{2m-n}=\Z$		 &&\\

		$\bullet$&	$\dfrac{a_{\lambda}^iu_{\lambda}^{m-i}}{u_{2\sigma}^{n/2}}\Big\rvert\dfrac{\bar a_{\lambda}^i\bar u_{\lambda}^{m-i}}{u_{\sigma}^{n}}\Big\rvert \zero$  &generates	& $H_{2m-n-2i}=\ev{\Z/4}$	 &for & $0<i<m$ \\
		
		$\bullet$&	$\dfrac{w_{2i+1}a_{\lambda}^m}{a_{\sigma}^{n-2i-1}}\Big\rvert \zero\Big\rvert \zero$ & generates &  $H_{-2i-1}=\ev{\Z/2}$  &for & $1\le i\le \frac n2-1$\\
		
		$\bullet$&	$\dfrac{2a_{\lambda}^m}{u_{2\sigma}^{n/2}}\Big\rvert \dfrac{\bar a_{\lambda}^m}{u_{\sigma}^n}\Big\rvert \zero$ &generates &$H_{-n}=Q$  &&
	\end{tabular}
\end{center}

\newpage
\subsection{\texorpdfstring{$\boldsymbol{H_*S^{m\lambda-n\sigma}}$ for odd $\boldsymbol n$}{Lambda minus Sigma nonoriented}}\label{Lambda minus sigma odd} If $m$ is nonzero,
	
	\begin{equation}
	H_*(S^{m\lambda-n\sigma})=
	\begin{cases}
	\hspace{-0.3em}
	\begin{tabular}{l p{0.5em} l p{1.5em} l p{1.5em} l}
	$\Z_{-}$			& if &  $*=2m-n\ge -1$		& &&&\\
	$\ev{\Z/2}\oplus \Z_{-}$		& if &$*=2m-n\le -3$	&&&&\\
	$\overline{\ev{\Z/2}}$		& if &$-1\le *<2m-n$ 		&and & $*$ is odd&&\\
	$\ev{\Z/2}$				& if & $2m-n<*\le -3$&and & $*$ is odd&&\\
	$\ev{\Z/2}\oplus \overline{\ev{\Z/2}}$				& if & $-n+2\le *<2m-n$&and & $*$ is odd &and & $*\le -3$\\
	$\ev{\Z/2}$				& if & $-n+1\le *<2m-n $&and & $*$ is even&&\\
	$Q$				& if & $*=-n$& and &$n\ge 3$&&\\
	$\overline{\ev{\Z/2}}$				& if & $*=-1$& and &$n=1$ & and &$m=1$\\
	\end{tabular}
	\end{cases}
	\end{equation}

	\begin{center}
		\begin{tabular}{p{0.2em} l p{3.8em} l p{1.5em} l}
			$\bullet$&	$\zero\Big\rvert u_{\sigma}^{-n}\bar u_{\lambda}^m\Big\rvert \bar u_{\sigma}^{-n}\bar{\bar u}_{\lambda}^m$		&		 generates &  $H_{2m-n}=\Z_{-}$		 &&\\

			$\bullet$&	$\zero\Big\rvert\dfrac{\bar a_{\lambda}^i\bar u_{\lambda}^{m-i}}{u_{\sigma}^{n}}\Big\rvert \zero$  &generates	& the $\overline{\ev{\Z/2}}$ in $H_{2m-n-2i}$	 &for & $0<i<m$ \\
			
			$\bullet$&	$\dfrac{w_{2i+1}a_{\lambda}^m}{a_{\sigma}^{n-2i-1}}\Big\rvert \zero\Big\rvert \zero$ & generates & the $\ev{\Z/2}$ in  $H_{-2i-1}$  &for & $1\le i\le \frac {n-3}2$\\
			
			$\bullet$&	$\dfrac{2a_{\lambda}^iu_{\lambda}^{m-i}}{a_{\sigma}u_{2\sigma}^{(n-1)/2}}\Big\rvert \zero\Big\rvert \zero$ &generates &$H_{2m-n-2i+1}=\ev{\Z/2}$ & for & $1\le i\le m$\\
			
			$\bullet$&	$w_na_{\lambda}^m \Big\rvert u_{\sigma}^{-n}\bar a_{\lambda}^m \Big\rvert \zero$ &generates &$H_{-n}=Q$ & for & $n\ge 3$\\
			
			$\bullet$&	$\zero\Big\rvert u_{\sigma}^{-1}\bar a_{\lambda}^m \Big\rvert \zero$ &generates &$H_{-1}=\overline{\ev{\Z/2}}$ & for & $n=m=1$
		\end{tabular}
	\end{center}

\newpage

\subsection{\texorpdfstring{$\boldsymbol{H_*S^{n\sigma-m\lambda}}$ for even $\boldsymbol n$}{Sigma minus Lambda oriented}}\label{Sigma minus Lambda even} If $n,m$ are both nonzero,

	\begin{equation}
H_*(S^{n\sigma-m\lambda})=
\begin{cases}
\hspace{-0.3em}
\begin{tabular}{l p{0.5em} l p{1.5em} l p{1.5em} l}
$Q^{\sharp}$			& if &  $*=n-3$		& and & $m\ge 2$&&\\
$\ev{\Z/4}$		& if & $n-2m<*<n-3$	& and & $*$ is odd&&\\
$\ev{\Z/2}$		& if &$0\le *\le n-4$&	and & $*$ is even& and &$*\neq n-2m$ 	\\
$L\oplus \ev{\Z/2}$				& if & $*=n-2m$& and &$n-2m\ge 0$&and &  $m\ge 2$\\
$L$				& if &  $*=n-2m$& and& $n-2m<0$ &and & $m\ge 2$\\
$L^{\sharp}$				& if &$*=n-2$  &and & $m=1$&&\\
\end{tabular}
\end{cases}
\end{equation}

\begin{center}
	\begin{tabular}{p{0.2em} l p{3.9em} l p{1.5em} l}
		$\bullet$&	$\dfrac{u_{2\sigma}^{n/2}s}{a_{\lambda}^{m-2}}\Big\rvert \dfrac{u_{\sigma}^n\bar s}{\bar a_{\lambda}^{m-2}}\Big\rvert \zero$		&		 generates &  $H_{n-3}=Q^{\sharp}$		 & for & $m\ge 2$\\

		$\bullet$&	$u_{2\sigma}^{n/2}\cdot \dfrac{s}{a_{\lambda}^{i-2}u_{\lambda}^{m-i}}\Big\rvert \dfrac{u_{\sigma}^n\bar s}{\bar a_{\lambda}^{i-2}\bar u_{\lambda}^{m-i}}\Big\rvert \zero$  &generates	& $H_{n-2m+2i-3}=\ev{\Z/4}$	 &for & $2\le i<m$ \\
		
		$\bullet$&	$\dfrac{a_{\sigma}^{n-2i}u_{2\sigma}^i}{a_{\lambda}^m}\Big\rvert \zero\Big\rvert \zero$ & generates & the $\ev{\Z/2}$ in  $H_{2i}$  &for &  $0\le i\le \frac{n-4}2$\\
		
		$\bullet$&	$\dfrac{4u_{2\sigma}^{n/2}}{u_{\lambda}^m}\Big\rvert \dfrac{2u_{\sigma}^n}{\bar u_{\lambda}^m}\Big\rvert  \dfrac{\bar u_{\sigma}^n}{\bar {\bar u}_{\lambda}^m}$ &generates &the  $L$ in $H_{n-2m}$  & for & $m\ge 2$\\
		
		$\bullet$&	$\dfrac{2u_{2\sigma}^{n/2}}{u_{\lambda}}\Big\rvert \dfrac{2u_{\sigma}^n}{\bar u_{\lambda}}\Big\rvert  \dfrac{\bar u_{\sigma}^n}{\bar {\bar u}_{\lambda}}$ &generates &$H_{n-2}=L^{\sharp}$  & for & $m=1$\\
	\end{tabular}
\end{center}
\phantom{}\bigbreak

\newpage

\subsection{\texorpdfstring{$\boldsymbol{H_*S^{n\sigma-m\lambda}}$ for odd $\boldsymbol n$}{Sigma minus Lambda nonoriented}}\label{Sigma minus Lambda odd} If $m$ is nonzero,

	\begin{equation}
H_*(S^{n\sigma-m\lambda})=
\begin{cases}
\hspace{-0.3em}
\begin{tabular}{l p{0.5em} l p{1.5em} l p{1.5em} l}
$Q^{\sharp}$			& if &  $*=n-3$		& and & $n\ge 3$& and &$ m\ge 2$\\
$\overline{\ev{\Z/2}}$		& if & $*=-2$	& and & $n=1$ & and &$m\ge 2$\\
$\ev{\Z/2}$		& if &$n-2m<*\le n-4$& and &$*$ is odd& & \\
$\ev{\Z/2}$		& if &$0\le *<n-2m$& and &$*$ is even& & \\
$\overline{\ev{\Z/2}}\oplus \ev{\Z/2}$				& if &  $n-2m<*\le n-5$ &and &  $*$ is even & and &$0\le *$ \\
$\overline{\ev{\Z/2}}$				& if &  $n-2m<*\le n-5$ &and &  $*$ is even& and &$*<0$ \\
$L_{-}$				& if &  $*=n-2m$& and& $m\ge 2$ &&\\
$\Z_{-}^{\flat}$				& if &$*=n-2$  &and & $m=1$&&
\end{tabular}
\end{cases}
\end{equation}

\begin{center}
	\begin{tabular}{p{0.2em} l p{3.8em} l p{1.5em} l}
		$\bullet$&	$\dfrac{a_{\sigma}^3u_{2\sigma}^{(n-3)/2}}{a_{\lambda}^m}\Big\rvert \dfrac{u_{\sigma}^n\bar s}{\bar a_{\lambda}^{m-2}}\Big\rvert \zero$		&		 generates &  $H_{n-3}=Q^{\sharp}$		 & for &  $n\ge 3$ and $m\ge 2$\\

		$\bullet$&	$\zero\Big\rvert \dfrac{u_{\sigma}\bar s}{\bar a_{\lambda}^{m-2}}\Big\rvert \zero$  &generates	& $H_{-2}=\overline{\ev{\Z/2}}$ 	 &for & $n=1$ and $m\ge 2$ \\
		
		$\bullet$&	$\dfrac{a_{\sigma}u_{2\sigma}^{(n-1)/2}s}{a_{\lambda}^{i-2}u_{\lambda}^{m-i}}\Big\rvert \zero\Big\rvert \zero$ & generates &  $H_{n-2m+2i-4}=\ev{\Z/2}$  &for &  $2\le i<m$\\
		
		$\bullet$&	$\dfrac{a_{\sigma}^{n-2i}u_{2\sigma}^{i}}{a_{\lambda}^m}\Big\rvert \zero\Big\rvert \zero$ & generates & the $\ev{\Z/2}$ in $H_{2i}$ &for & $0\le i< \frac{n-5}2$\\
		
		$\bullet$&	$\zero\Big\rvert \dfrac{u_{\sigma}^n\bar s}{\bar a_{\lambda}^{i-2}\bar u_{\lambda}^{m-i}}\Big\rvert \zero$  &generates	& the $\overline{\ev{\Z/2}}$ in $H_{n-2m+2i-3}$	 &for & $2\le i<m$  \\
				
		$\bullet$&	$\dfrac{a_{\sigma}u_{2\sigma}^{(n-1)/2}s}{u_{\lambda}^{m-2}}\Big\rvert \dfrac{2u_{\sigma}^n}{\bar u_{\lambda}^m}\Big\rvert  \dfrac{\bar u_{\sigma}^n}{\bar {\bar u}_{\lambda}^m}$ &generates &$H_{n-2m}=L_{-}$ & for & $m\ge 2$\\
		
		$\bullet$&	$\zero\Big\rvert \dfrac{2u_{\sigma}^n}{\bar u_{\lambda}}\Big\rvert  \dfrac{\bar u_{\sigma}^n}{\bar {\bar u}_{\lambda}}$ &generates &$H_{n-2}=\Z_{-}^{\flat}$  & for & $m=1$\\
	\end{tabular}
\end{center}

\newpage
\section{\texorpdfstring{Computing the $RO(G)$ graded homology}{Computing the RO(G) graded homology}}\label{Framework}
We explain the theoretical framework behind our computation of the $RO(G)$ homology of a point for a finite group $G$. The results here are classical, but including them allows us to expound on the chain complexes and spectral sequences that we'll be using in the following sections. \medbreak

First, let us recall that the shift of a $G$-Mackey functor $M$ at a finite $G$-set $T$ is the Mackey functor $M_T$ specified on orbits as 
\begin{equation}
M_T(G/H)=M(T\times G/H)
\end{equation}

Our goal is to compute 
\begin{equation}
H_*(S^V;M)
\end{equation}
where $V$ is a virtual representation of $G$; in particular $S^V$ is a spectrum, not a space (not even necessarily a suspension spectrum). For a general finite $G$-spectrum $X$,
\begin{equation}
H_*(X;M)
\end{equation}
is computed using an equivariant cell decomposition of $X$. This is a sequence of $G$-spectra $X_p$ interpolating $X_0=S$ and $X_n=X$ through cofiber sequences
\begin{equation}
X_{p-1}\to X_p\to T_{p+}\wedge S^p
\end{equation}
where the $T_p$ are finite $G$-sets. Given such a decomposition, we have an Atiyah-Hirzebruch spectral sequence of Mackey functors
\begin{equation}
E^1_{p,q}=H_q(T_{p+},M)\implies H_{p+q}(X,M)
\end{equation}
But by the definition of equivariant homology, $H_q(T_{p+},M)=M_{T_p}$ concentrated in degree $q=0$, hence the $E_1$ page is actually a chain complex with
\begin{equation}
C_p(X;M)=M_{T_p}
\end{equation}
The boundary maps are induced from the geometric boundary maps 
\begin{equation}
T_{p+}\wedge S^p\to \Sigma X_{p-1}\to \Sigma(T_{(p-1)+}\wedge S^{p-1})
\end{equation}
 in the following way: Smash the composite above with $S^{-p}$ to get a $G$-map $T_{p}\to T_{p-1}$ and then use the induced transfer map $M_{T_p}\to M_{T_{p-1}}$ that is specified as:
\begin{equation}
M_{T_p}(G/H)=M(G/H\times T_p)\xrightarrow{\Tr}M(G/H\times T_{p-1})=M_{T_{p-1}}(G/H)
\end{equation}
This is the algebraic boundary map $C_p\to C_{p-1}$.

The homology of the chain complex $C_*(X;M)$ is $H_*(X;M)$. We can do the same for cohomology and get the dual cochain complex $C^*(X;M)$ (using the induced restriction maps $C^p=M_{T_p}\to  C^{p+1}=M_{T_{p+1}}$).\medbreak

So if we have an equivariant cell decomposition of $X$ the problem of computing $H_*(X;M)$ is reduced to algebra. If $V$ is an actual (as opposed to virtual) representation, we might be able to find this decomposition of $S^V$ from the geometry of the space, and that's what we do in section \ref{ProofsPureHomology}. If $V=-W$ where $W$ is an actual representation then we can use Spanier-Whitehead duality and be reduced to the case already considered:
\begin{equation}
H_*(S^V;M)=H^{-*}(S^W;M)
\end{equation}
However, we can't perform this trick if the virtual representation is $V-W$. Instead, we use that $S^{V-W}=S^V\wedge S^{-W}$ and smash the cell decompositions for $S^V,S^W$ together to get one for $S^{V-W}$. In general, given cell decompositions
\begin{gather}
X_{n-1}\to X_n\to T_{n+}\wedge S^n\\
Y_{n-1}\to Y_n\to T_{n+}'\wedge S^n
\end{gather}
we get a cell decomposition of $X\wedge Y$ by
\begin{equation}
ho\colimit_{k+l=n}X_k\wedge Y_l\to ho\colimit_{k+l=n+1}X_k\wedge Y_l\to \Big(\coprod_{k+l=n+1}T_k\times T'_l\Big)_+\wedge S^{k+l}
\end{equation}
(this reduces to a fact in symmetric monoidal triangulated categories proven in \cite{May01}). Therefore,
\begin{gather}
C_*(X\wedge Y;M\boxtimes N)=(M\boxtimes N)_{\coprod_{k+l=*}T_k\times T'_l}=\oplus_{k+l=*}(M\boxtimes N)_{T_k\times T'_l}=\\
=\oplus_{k+l=*}M_{T_k}\boxtimes N_{T'_l}=C_*(X;M)\boxtimes C_*(Y;N)
\end{gather}
The boundary maps match up as well.\medbreak

In our case, we take $M=N=\Z$ and then $H_*(S^{V-W};\Z)$ is computed as the homology of
\begin{equation}
C_*(S^V;\Z)\boxtimes C^{-*}(S^W;\Z)
\end{equation}
Unfortunately, even if the $C_*S^V$ and $C^{-*}S^W$ are by themselves very small and easy to compute with (like in section \ref{ProofsPureHomology}), the box product can be extremely large very easily. So instead of a direct computation, we use algebraic spectral sequences converging to its homology.\medbreak

In general, for any tensor product of chain complexes $C\otimes D$ in a sufficiently good symmetric monoidal abelian category (like that of Mackey functors), we have three spectral sequences converging to $H_*(C\otimes D)$. If we filter the double complex underlying the tensor product either horizontally or vertically, we get two spectral sequences with $E_2$ terms
\begin{gather}
E_2=H_*(C;H_*D)\implies H_*(C\otimes D)\\
E_2=H_*(D;H_*C)\implies H_*(C\otimes D)
\end{gather}
Using Cartan-Eilenberg resolutions we obtain a Kunneth spectral sequence
\begin{gather}
E_2=\Tor^{*,*}(H_*C,H_*D)\implies H_*(C\otimes D)
\end{gather}
We refer the reader to \cite{Wei94} and \cite{Rot09} for details on their constructions but we shall not need them.

In our case of $C_*(S^V)\boxtimes C^{-*}(S^W)$ the spectral sequences take the form
\begin{gather}
E^2_{p,q}=H_p(S^V,H^{-q}S^W)\implies H_{p+q}S^{V-W}\\
E_2^{p,q}=H^p(S^W,H_{-q}S^V)\implies H_{p+q}S^{V-W}\\
E^2_{p,q}=\Tor^{p,q}(H_pS^V,H_qS^{-W})\implies H_{p+q}S^{V-W}
\end{gather}
These are all spectral sequences of $\underline{\Z}$-module and the final one uses the $\Tor$ in the symmetric monoidal category of  $\underline{\Z}$-modules.  

Finally, we remark that our three spectral sequences can also be obtained topologically. The first two are the Atiyah-Hirzebruch spectral sequences for the homology theory $H_V$ and the final one is the topological Kunneth spectral sequence.  

%\newcounter{proofpart}
%\stepcounter{proofpart}

\section{\texorpdfstring{Proofs for $H_*(S^{n\sigma+m\lambda})$}{Proofs for nsigma+mlambda}}\label{ProofsPureHomology}

The results of this section also appear in \cite{HHR17}, but we have chosen to include them for the sake of completeness. As always, $n,m\ge 0$.\medbreak

We can obtain equivariant cell decompositions for $S^{n\sigma+m\lambda}$ as follows. View it as the compactification of the disc $D(\R^n\times \R^{2m})$ and include the $C_4$ subspace $X_{n+2m-1}$ where either one of the final two coordinates is $0$. The quotient space is 
\begin{equation}
S^{n\sigma+m\lambda}/X_{n+2m-1}=C_{4+}\wedge S^{n+2m}
\end{equation}
This is the wedge of four $S^{n+2m}$'s that correspond to the signs of the last two coordinates in $D(\R^n\times \R^{2m})$, and we use $(x_{\pm},y_{\pm})$ to represent them.

The space $X_{n+2m-1}$ includes $S^{n\sigma+(m-1)\lambda}$ as the $C_4$ subspace where both the last two coordinates are $0$. The quotient is $C_{4+}\wedge S^{n+2m-1}$ and we use $(x_{\pm},0)$ and $(0,y_{\pm})$ to represent the four spheres.

We continue like this until we reach $S^{n\sigma}$ and we then include $S^{(n-1)\sigma}$ (last coordinate $0$) with quotient $(C_4/C_2)_+\wedge S^n$; these are two spheres represented by $x_+,x_-$. Eventually we will reach $S^0$ represented by "$1$". \medbreak

We write $\Z_{C_4},\Z_{C_2}$ for the shifts of $\underline{\Z}$ at the orbits $C_4/C_4$ and $C_4/C_2$; these are also the fixed point Mackey functors of $\Z[C_4], \Z[C_2]$ respectively.\medbreak

In this notation, the chain complex for $S^{n\sigma+m\lambda}$ is
\begin{equation}
0\to \Z_{C_4}\to \Z_{C_4}\to \cdots\to \Z_{C_4}\to \Z_{C_2}\to \Z_{C_2}\to \cdots\to \Z_{C_2}\to \Z\to 0
\end{equation}
The $\Z_{C_2}$'s are generated by $x_{\pm}$ while the $\Z_{C_4}$'s are interchangeably generated by $(x_{\pm},y_{\pm})$ and $(x_{\pm},0),(0,y_{\pm})$.

The differentials up to $d_n$ are $d_1x_+=1$, $d_{2k}x_+=x_+-x_-$ and $d_{2k+1}=x_++x_-$.

The differentials from $d_{n+1}$ to $d_{2m+n}$ depend on whether our sphere is $C_4$-oriented or not. If $n$ is even (oriented),  
\begin{gather}
d(x_{+},y_+)=(x_+,0)-(0,y_+)\\
d(x_+,0)=\sum^4 g^i(x_+,y_+)
\end{gather}

If $n$ is odd (non oriented), 
\begin{gather}
d(x_{+},y_+)=(x_+,0)+(0,y_+)\\
d(x_+,0)=\sum^4 (-1)^ig^i(x_+,y_+)
\end{gather}
The differential $d_{n+1}:\Z_{C_4}\to \Z_{C_2}$ is $d(x_+,0)=x_++x_-$ if $n$ is even and $d(x_+,0)=x_+-x_-$ if $n$ is odd. Finally, if $n=0$ the differential $d_1:\Z_{C_4}\to \Z$ is $d_1(x_+,0)=1$.

These differentials can be computed geometrically, or inferred from the observation that the homology of the bottom level of $C_*S^{n\sigma+m\lambda}$ is the nonequivariant homology of $S^{n\sigma+m\lambda}$ i.e. $S^{n+2m}$, so the bottom level of $C_*S^{n\sigma+m\lambda}$ must be exact apart from the highest degree.\medbreak

We readily compute the homology of $C_*S^{n\sigma+m\lambda}$ to be as follows:\medbreak

If $n$ is even,
\begin{gather}
H_*(S^{n\sigma+m\lambda})=\begin{cases}
\Z&\textup{if $*=n+2m$}\\
\ev{\Z/4}&\textup{if $n\le *<n+2m$ and $*$ is even}\\
\ev{\Z/2}&\textup{if $0\le *<n$ and $*$ is even}
\end{cases}
\end{gather}

If $n$ is odd,
\begin{gather}
H_*(S^{n\sigma+m\lambda})=\begin{cases}
\Z_{-}&\textup{if $*=n+2m$}\\
\overline{\ev{\Z/2}}&\textup{if $n\le *<n+2m$ and $*$ is odd}\\
\ev{\Z/2}&\textup{if $0\le *<n+2m$, $*$ is even}
\end{cases}
\end{gather}

We now describe the multiplicative generators of  $H_*(S^{n\sigma+m\lambda})$. Recall from section \ref{Generators} that we have the Euler and orientation classes generating the Mackey functors:
\begin{gather}
a_{\sigma}|\zero|\zero\in H_0(S^{\sigma})=\ev{\Z/2}\\
a_{\lambda}|\bar a_{\lambda}|\zero\in H_0(S^{\lambda})=\ev{\Z/4}\\
u_{2\sigma}|u_{\sigma}^2|\bar u_{\sigma}^2\in H_2(S^{2\sigma})=\Z\\
u_{\lambda}|\bar u_{\lambda}|\bar{\bar{u}}_{\lambda}\in H_2(S^{\lambda})=\Z\\
\zero|u_{\sigma}|\bar u_{\sigma}\in H_1(S^{\sigma})=\Z_{-}
\end{gather}
and satisfying the Gold Relation:
\begin{gather}a_{\sigma}^2u_{\lambda}=2a_{\lambda}u_{2\sigma}\end{gather}	

%Note that for even $n$, if we multiply the $\ev{\Z/2}$ generator by $u_{\lambda}$ we'll get either the $\ev{\Z/2}$ generator (if $i\le n/2-2$) or twice the $\ev{\Z/4}$ generator (if $i=n/2-1$) in a higher sphere, by the Gold Relation. This is why we don't have $a_{\sigma},a_{\lambda},u_{\lambda}$ all together in that case.\medbreak

These classes multiplicatively generate all of $H_*(S^{n\sigma+m\lambda})$ and the only primary relation is the Gold. This claim follows easily from the following observations: Multiplication by $a_{\sigma}:S^0\to S^{\sigma}$ induces the chain map
\begin{equation}
C_*S^{n\sigma+m\lambda}\to C_*S^{(n+1)\sigma+m\lambda}
\end{equation}
that is $\Z_{C_2}\xrightarrow{1}\Z_{C_2}$ for $*\le n$, $\Z_{C_4}\xrightarrow{1}\Z_{C_4}$ for $n+1\le *\le n+2m$ and the map $\Z_{C_4}\to \Z_{C_2}$ given by the canonical projection $\Z[C_4]\to \Z[C_4/C_2]$ %($x\mapsto x$)
at bottom level. For $m=0$, $a_{\sigma}$ is the canonical inclusion $C_*S^{n\sigma}\to C_*S^{(n+1)\sigma}$. Similarly, multiplication by $a_{\lambda}$ induces the canonical inclusion $C_*S^{n\sigma+m\lambda}\to C_*S^{(n+1)\sigma+m\lambda}$.

From these observations it follows that multiplication by $a_{\sigma},a_{\lambda}$ is an isomorphism in certain dimensions which is enough to prove the multiplicative generation of $H_*(S^{n\sigma+m\lambda})$ by Euler and orientation classes.

%\stepcounter{proofpart}
\section{\texorpdfstring{Proofs for  $H_*(S^{-n\sigma-m\lambda})$}{Proofs for -nsigma-mlambda}}\label{ProofsPureCohomology}

To get $H_*(S^{-n\sigma-m\lambda})$ we use Spanier-Whitehead Duality:
\begin{equation}
H_*(S^{-n\sigma-m\lambda})=H^{-*}(S^{n\sigma+m\lambda})
\end{equation}
and cohomology is computed by the dual chain complex $C^*(S^{n\sigma+m\lambda})$. The results:

If $n$ is even,

\begin{gather}
H_*(S^{-n\sigma-m\lambda})=\begin{cases}
L&\textup{if $*=-n-2m$ and $m\neq 0$}\\
p^*L&\textup{if $*=-n-2m$ and $m=0$}\\
\ev{\Z/4}&\textup{if $-n-2m<*<-n-1$ and $*$ is odd}\\
\ev{\Z/2}&\textup{if $-n-1\le *<-1$ and $*$ is odd and $m\neq 0$}\\
\ev{\Z/2}&\textup{if $-n+1\le *<-1$ and $*$ is odd and $m=0$}
\end{cases}
\end{gather}

If $n$ is odd,

\begin{gather}
H_*(S^{-n\sigma-m\lambda})=\begin{cases}
L_{-}&\textup{if $*=-n-2m<-1$ and $m\neq 0$}\\
p^*L_{-}&\textup{if $*=-n<-1$ and $m=0$}\\
\Z_{-}&\textup{if $*=-1$ and $n=1$ and $m=0$}\\
\overline{\ev{\Z/2}}&\textup{if $-n-2m<*<-n-1$ and $*$ is even}\\
\ev{\Z/2}&\textup{if $-n-2m<*<-1$, $*$ is odd}
\end{cases}
\end{gather}

%The $n=0$ computation is consistent with $H^*_{C_4}(S^{m\lambda})=H^{*-1}_{C_4}(S(m\lambda)_+)$ for $*\ge 2$ and $S(m\lambda)_+/C_4=L(4;1,...,1)_+$ the generalized Lens space. The \emph{integral} cohomology is $\Z/4$ in every even degree from $2$ to $2m-2$. \medbreak

We shall now find the multiplicative generators of $H_*(S^{-n\sigma-m\lambda})$.

First note that $u_{\sigma}\in H_2S^{\lambda}=\Z\{x_+-x_-\}$ pairs with the generator of $H_{-2}S^{-\lambda}=\Z_{-}\{x_+^*\}$ to $x_+^*(x_+-x_-)=1$ hence it's invertible. Similarly $\bar{\bar{u}}_{\lambda}$ is invertible.

The transfers on the products of $u_{\sigma}^{-1},\bar{\bar{u}}_{\lambda}^{-1}$ generate the $L,L_{-}$'s since these elements do so on the bottom level (and transfers are surjective for $L,L_{-}$).  We can compute these transfers through the Frobenius relations; for example,
\begin{equation}
\Tr_1^4(\bar u_{\sigma}^{-2n}\bar{\bar{u}}_{\lambda}^{-m})=4/(u_{2\sigma}^{n}u_{\lambda}^{m})
\end{equation}
We have:
\begin{itemize}
	\item $4/(u_{2\sigma}^{n}u_{\lambda}^{m})|2/(u_{\sigma}^{2n}\bar u_{\lambda}^m)|\bar u_{\sigma}^{-2n}\bar{\bar{u}}_{\lambda}^{-m}$ generates $L$ for $m>0$
	\item $2/u_{2\sigma}^{n}|u_{\sigma}^{-2n}|\bar u_{\sigma}^{-2n}$ generates $p^*L$.
\end{itemize}

For odd $n\ge 3$  let
\begin{gather}
w_n=\Tr_2^4(u_{\sigma}^{-n})=w_3/u_{2\sigma}^{(n-3)/2}
\end{gather}	
generating the top level of  $H_{-n}S^{-n\sigma}=p^*L_{-}$. For $n=1$ this transfer is $0$ so we don't define a $w_1$. Next, for odd $n\ge 1$ and any $m\ge 1$ let
\begin{equation}
x_{n,m}=\Tr_1^4(\bar u_{\sigma}^{-n}\bar{\bar{u}}_{\lambda}^{-m})=x_{1,1}/(u_{2\sigma}^{(n-1)/2}u_{\lambda}^{m-1})
\end{equation}
generating the top level of $H_{-n-2m}S^{-n\sigma-m\lambda}=L_{-}$. We have
\begin{itemize}
\item $x_{n,m}|2/(u_{\sigma}^n\bar{u}_{\lambda}^m)|\bar u_{\sigma}^{-n}\bar{\bar{u}}_{\lambda}^{-m}$ generates $L_{-}$.
\item $w_n|u_{\sigma}^{-n}|\bar u_{\sigma}^{-n}$ generates $p^*L_{-}$.
\end{itemize}

The $w_n$ are infinitely divisible by $a_{\sigma}$ since $a_{\sigma}:S^{n\sigma}\to S^{(n+1)\sigma}$ is the inclusion in chains $C_*S^{n\sigma}\to C_*S^{(n+1)\sigma}$ hence projection in cochains $C^*S^{(n+1)\sigma}\to C^*S^{n\sigma}$ which is a quasi-isomorphism in top level for odd $n$. 

\begin{itemize}
	\item $w_n/a_{\sigma}^i|\zero|\zero$ generates $\ev{\Z/2}$
\end{itemize}	

The $w_n/a_{\sigma}^i$ are infinitely divisible by $a_{\lambda}$ since $a_{\lambda}:S^{n\sigma+m\lambda}\to S^{n\sigma+(m+1)\lambda}$ gives the projection of cochain complexes $C^*S^{n\sigma+(m+1)\lambda}\to C^*S^{n\sigma+m\lambda}$ which is a quasi-isomorphism in top level.
%. Therefore it is an isomorphism in $C^i$ for $i\le n+2m$, and becomes an iso on $H_{>n-2m}$ and inclusion on $H_{-n-2m}$. This is an isomorphism on top level.\medbreak

Now note that the $x_{n,m}$ are also infinitely $a_{\sigma}$ divisible. This is because 
\begin{equation}
a_{\sigma}:C^*S^{(n+1)\sigma+m\lambda}\to C^*S^{n\sigma+m\lambda}
\end{equation}
is identity for $*\le n$ or $n+2\le *\le n+2m$, while for $*=n+1$ it's the map $\Z_{C_2}\to \Z_{C_4}$ dual to $\Z[C_4]\to \Z[C_4/C_2]$. This is an isomorphism in top level.

The $w_n/(a_{\sigma}^ia_{\lambda}^j)$ and $x_{n,m}/a_{\sigma}^i$ generate the top levels of $\ev{\Z/2}$'s, so if they occur at the same dimensions they must be equal. Thus
\begin{equation}x_{1,1}/a_{\sigma}^{2}=w_3/a_{\lambda}
\end{equation}
and in general,
\begin{equation}
x_{n,m}/a_{\sigma}^i=x_{n,1}/(a_{\sigma}^{i}u_{\lambda}^{m-1})=w_{n+2}/(a_{\sigma}^{i-2}a_{\lambda}u_{\lambda}^{m-1})
\end{equation}
for odd $n$ and $m\ge 1$ and $i\ge 2$.  

So the $x_{n,1}/a_{\sigma}^{i}$ are infinitely $a_{\lambda}$ divisible for $i\ge 2$ and odd $n$.
\begin{itemize}
	\item $x_{2i+1,1}/(a_{\sigma}^{n-2i-1}a_{\lambda}^{m-1})|\zero|\zero$ generates $\ev{\Z/2}$.
\end{itemize}

Let $s\in H_{-3}(S^{-2\lambda})=\ev{\Z/4}$ be the generator (unique up to a sign). It is infinitely divisible by $u_{\lambda}$: Indeed, we have the commutative diagram
\begin{center}\begin{tikzcd}
H_{-3-2i}(C_4/C_4)=\Z/4\ar[r,"u_{\lambda}^i"]\ar[d,"\Res^4_2"]&\Z/4=H_{-3}(C_4/C_4)\ar[d,"\Res^4_2"]\\
H_{-3-2i}(C_4/C_2)=\Z/2\ar[r,"\bar u_{\lambda}^i"]&\Z/2=H_{-3}(C_4/C_2)
\end{tikzcd}\end{center}
where the right column is generated by $s,\bar s$ and the bottom horizontal map is an isomorphism by the $\pi^{C_2}_{\bigstar}(H\underline{\Z})$ computation in \cite{HHR17} (because $C_2$ only has one nontrivial irreducible real representation this computation is significantly shorter and easier than the $C_4$ case, but we will not reproduce it here).

The element $s/u_{\lambda}^i$ is infinitely $a_{\lambda}$ divisible (by either the $C_2$ restriction argument or the usual inclusion of chains argument) and $s/(a_{\lambda}^{i}u_{\lambda}^{j})$ generates the top level of a $\ev{\Z/4}$. To get the remaining $\ev{\Z/4}$ which appear as in the homology of $S^{-n\sigma-m\lambda}$ for even $n$ we can further multiply by $u_{\sigma}^{-n}$ in the middle level, which gives $s/(u_{2\sigma}^{n/2}a_{\lambda}^{i}u_{\lambda}^{j})$ as the top level generator. We conclude:
\begin{itemize}
	\item $s/(u_{2\sigma}^{n/2}a_{\lambda}^{i}u_{\lambda}^{j})|\bar s/(u_{\sigma}^{n}\bar a_{\lambda}^{i}\bar u_{\lambda}^{j})|\zero$ generates $\ev{\Z/4}$ for $i,j\ge 0$ and even $n$.
	\item $\zero|\bar s/(u_{\sigma}^{n}\bar a_{\lambda}^{i}\bar u_{\lambda}^{j})|\zero$ generates $\overline{\ev{\Z/2}}$ for $i,j\ge 0$ and odd $n$.
\end{itemize}

%\stepcounter{proofpart}
\section{Preparation for \texorpdfstring{$H_*(S^{m\lambda-n\sigma})$}{the case m*lambda-n*sigma}}\label{ProofsE2page1}

In this section and the next we will compute $H_*(S^{m\lambda-n\sigma})$ for $n,m\ge 1$.

Recall from section \ref{Framework} that there are three algebraic spectral sequences we can use for this computation. This section is devoted to determining the Mackey functors in their $E_2$ pages; we will figure out the differentials and extensions in the next section.

We omit the multiplicative presentation of certain generators in the $E_2$ page of some of our spectral sequences. These generators either don't survive the spectral sequence, or if they do then we can figure out their multiplicative presentation in the $E_{\infty}$ page by comparison with the other spectral sequences.

\subsection{The Homological  Spectral Sequence}

The HSS  is
\begin{equation}
E^2_{p,q}=H_p(S^{m\lambda};H_qS^{-n\sigma})\implies H_{p+q}(S^{m\lambda-n\sigma})
\end{equation}
We compute the $E_2$ page (sans differentials) directly with the chain complex 
\begin{equation}
C_*(S^{m\lambda};H_qS^{-n\sigma})=C_*(S^{m\lambda})\boxtimes H_q(S^{-n\sigma})
\end{equation}
and get the following results. First,
\begin{gather}
H_*(S^{m\lambda};H_{-1}S^{-\sigma})=H_*(S^{m\lambda};\Z_{-})=\begin{cases}
\Z_{-}&\textup{if $*=2m$}\\
\overline{\ev{\Z/2}}&\textup{if $0\le *\le 2m-2$ and $*$ is even}\\
\ev{\Z/2}&\textup{if $1\le *\le 2m-1$ and $*$ is odd} \end{cases}
\end{gather}
The multiplicative generators:
\begin{itemize}
	\item $\zero|u_{\sigma}^{-1}\bar u_{\lambda}^m|\bar u_{\sigma}^{-1}\bar{\bar u}_{\lambda}^m$ generates $\Z_{-}$
	\item $\zero|u_{\sigma}^{-1}\bar a_{\lambda}^i\bar u_{\lambda}^{m-i}|\zero$ generates $\overline{\ev{\Z/2}}$ for $0<i\le m$
\end{itemize}
\phantom{}\smallbreak
Second, for even $n\ge 2$,
\begin{equation}
H_*(S^{m\lambda};H_{-n}S^{-n\sigma})=H_{*}(S^{m\lambda};p^*L)=\begin{cases}
\Z&\textup{if $*=2m$}\\
\ev{\Z/4}&\textup{if $2\le *<2m$ and $*$ is even}\\
Q&\textup{if $*=0$}
\end{cases}
\end{equation}
The multiplicative generators:
\begin{itemize}
	\item $u_{\lambda}^m/u_{2\sigma}^{n/2}|u_{\sigma}^{-n}\bar u_{\lambda}^m|\bar u_{\sigma}^{-n}\bar{\bar{u}}_{\lambda}^m$ generates $\Z$
	\item $(a_{\lambda}^iu_{\lambda}^{m-i})/u_{2\sigma}^{n/2}|\bar u_{\sigma}^{-n}\bar a_{\lambda}^i\bar u_{\lambda}^{m-i}|\zero$ generates $\ev{\Z/4}$ for $0<i<m$
	\item $(2a_{\lambda}^m)/u_{2\sigma}^{n/2}|u_{\sigma}^{-n}\bar a_{\lambda}^m|\zero$ generates $Q$.
\end{itemize}
\phantom{}\smallbreak
Third, for odd $n\ge 3$,
\begin{equation}
H_*(S^{m\lambda};H_{-n}S^{-n\sigma})=H_*(S^{m\lambda};p^*L_{-})=\begin{cases}
\Z_{-}&\textup{if $*=2m$}\\
\ev{\Z/2}&\textup{if $1\le *\le 2m-1$ and $*$ is odd}\\
\overline{\ev{\Z/2}}&\textup{if $2\le *\le 2m-2$ and $*$ is even}\\
Q&\textup{if $*=0$}
\end{cases}
\end{equation}
The multiplicative generators:
\begin{itemize}
	\item $\zero|u_{\sigma}^{-n}\bar u_{\lambda}^m|\bar u_{\sigma}^{-n}\bar{\bar{u}}_{\lambda}^m$ generates $\Z_{-}$
	\item $\zero|\bar u_{\sigma}^{-n}\bar a_{\lambda}^i\bar u_{\lambda}^{m-i}|\zero$ generates $\overline{\ev{\Z/2}}$ for $0<i<m$
	\item $w_na_{\lambda}^m|u_{\sigma}^{-n}\bar a_{\lambda}^m|\zero$ generates $Q$.
\end{itemize}
\phantom{}\medbreak

Finally,  
\begin{equation}
H_*(S^{m\lambda};H_{-2i-1}S^{-n\sigma})=H_*(S^{m\lambda};\ev{\Z/2})=\ev{\Z/2}
\end{equation}
concentrated in degree $*=0$ and generated by $(w_{2i+1}a_{\lambda}^m)/a_{\sigma}^{n-2i+1}$ for $1\le i<(n-1)/2$.\bigbreak

\subsection{The Cohomological  Spectral Sequence}

The CSS  is
\begin{equation}
E_2^{p,q}=H_{-p}(S^{-n\sigma};H_{-q}S^{m\lambda})\implies H_{-p-q}(S^{m\lambda-n\sigma})
\end{equation}
The Mackey functors in the $E_2$ page are computed in a similar fashion as the HSS, and are as follows: First, $H_{-*}(S^{-n\sigma};H_{2m}S^{m\lambda})=H_{-*}(S^{-n\sigma};\Z)$ was computed in the previous section. Next,
\begin{gather}
H_{-*}(S^{-n\sigma};H_iS^{m\lambda})=H_{-*}(S^{-n\sigma};\ev{\Z/4})=\begin{cases}
Q&\textup{if $*=n\ge 2$}\\
\ev{\Z/2}&\textup{if $0\le *<n$ and $*\neq 1$}\\
\overline{\ev{\Z/2}}&\textup{if $*=n=1$}
\end{cases}
\end{gather}
for even $i$ is even and $0\le i<2m$. The multiplicative generators:

\begin{itemize}  \setlength\itemsep{0.5em}
	
	\item   $(2a_{\lambda}^{m-i}u_{\lambda}^i)/u_{\sigma}^{n/2}|u_{\sigma}^{-n}\bar a_{\lambda}^{m-i}\bar u_{\lambda}^i|\zero$ generates $Q$ for even $n\ge 2$
	\item $w_na_{\lambda}^{m-i}u_{\lambda}^i|u_{\sigma}^{-n}\bar a_{\lambda}^{m-i}\bar u_{\lambda}^i|\zero$ generates $Q$ for odd $n\ge 3$.
	\item $(w_{2j+1}a_{\lambda}^{m-i}u_{\lambda}^i)/a_{\sigma}^{n-2j-1}|\zero|\zero$ generates $\ev{\Z/2}$ in odd degrees.
	\item $(2a_{\lambda}^{m-i}u_{\lambda}^i)/a_{\sigma}^n|\zero|\zero$ generates $\ev{\Z/2}$ at degree $0$. % This is because  $a_{\sigma}^n:C^*(S^{n\sigma};\ev{\Z/4})\to C^*(S^0;\ev{\Z/4})$  is the projection map which is an isomorphism for $*=0$.
	\item $\zero|u_{\sigma}^{-1}\bar a_{\lambda}^{m-i}\bar u_{\lambda}^i|\zero$ generates $\overline{\ev{\Z/2}}$
\end{itemize}

%I once wrote that  $u_{2\sigma}^{-j}a_{\sigma}^{-n+2j}(2a_{\lambda}^{m-i})u_{\lambda}^i$ generates for every even degree $0\le *<n$ ($0\le j<n/2$). The problem with the proof is that $u_{2\sigma}$ is not a map of chains and can only be seen in homology, say using Kunneth. The elements do exist in the mixed homology for that reason, but perhaps not in the twisted homology (not sure).

\subsection{The Kunneth Spectral Sequence}

The KSS is
\begin{equation}
E^2_{p,q}=\Tor^{p,q+1}_{\underline{\Z}}(H_*S^{m\lambda-(n-1)\sigma},\Z_{-})\implies H_{p+q}(S^{m\lambda-n\sigma})
\end{equation}
This is the Kunneth spectral sequence for $S^{m\lambda-(n-1)\sigma}\wedge S^{-\sigma}$ and we have used that $H_*S^{-\sigma}=\Z_{-}$ concentrated in degree $*=-1$. A free resolution of $\Z_{-}$ is
\begin{equation}
0\to \Z\xrightarrow{\Delta}\Z_{C_2}\xrightarrow{\nabla} \Z_{-}\to 0
\end{equation}
We list only the nonzero $\Tor$'s below, dropping the subscript $\underline{\Z}$ from the notation:

\begin{itemize}  \setlength\itemsep{0.5em}
	\item $\Tor^0(\ev{\Z/4},\Z_{-})=\overline{\ev{\Z/2}}$ and $\Tor^1(\ev{\Z/4},\Z_{-})=\ev{\Z/2}$
	\item $\Tor^0(\overline{\ev{\Z/2}},\Z_{-})=Q$
	%	($\overline{\ev{\Z/2}}_{C_2}=\ev{\Z/4}_{C_2}$)
	\item $\Tor^0(Q,\Z_{-})=Q$ and $\Tor^1(Q,\Z_{-})=\ev{\Z/2}$ %($Q_{C_2}=\ev{\Z/4}_{C_2}$)
	\item $\Tor^1(\ev{\Z/2},\Z_{-})=\ev{\Z/2}$  %($\ev{\Z/2}_{C_2}=0$)
	\item $\Tor^0(\Z_{-},\Z_{-})=p^*L$ %$\ev{\Z_{-}}_{C_2}=\Z_{C_2}$ so 
	\item $\Tor^0(p^*L,\Z_{-})=p^*L_{-}$ % $(p^*L)_{C_2}=\Z_{C_2}$ % 
	\item $\Tor^0(p^*L_{-},\Z_{-})=p^*L$ and $\Tor^1(p^*L_{-},\Z_{-})=\ev{\Z/2}$ %$(p^*L_{-})_{C_2}=\Z_{C_2}$ but on a new generator y with y=x, gy=-gx (!)
	\item $\Tor^0(L,Z_{-})=L_{-}$
	\item $\Tor^0(L_{-},\Z_{-})=L$ and $\Tor^1(L_{-},\Z_{-})=\ev{\Z/2}$
\end{itemize}

Note that $\Tor^*(\Z_{-},-)$ vanishes above $*=1$ by the resolution of $\Z_{-}$. This means that the KSS is concentrated in the first two columns, hence \emph{always collapses} for dimensional reasons. That said, we still need the HSS to solve the extension problems in the KSS and get the multiplicative generators. We only make a single use of the CSS and that's for $n=1$.

\section{Proofs for \texorpdfstring{$H_*(S^{m\lambda-n\sigma})$}{The case m*lambda-n*sigma}}\label{Proofssigma}

We first compute $H_*(S^{m\lambda-n\sigma})$ separately for $n=1$ through $n=5$ using the spectral sequences explained in the preceding section. 
%With the exception of $n=1$, which is quickly dealt with by the CSS, all other cases follow by comparing the HSS with the KSS, using one to determine the differentials and extensions of the other. 
The KSS used to compute $S^{m\lambda-n\sigma}$ works by feeding it the answer of the computation for $n-1$, so we perform our calculations in order of increasing $n$. The general $n$ case follows by induction, exhibiting very similar behavior to the $n=4$ and $n=5$ cases, depending on the parity of $n$.

\subsection{ The case $n=1$}
The CSS for $S^{m\lambda-\sigma}$ is concentrated in the first two rows and even columns, hence collapses with no extensions to 
\begin{equation}
H_*(S^{m\lambda-\sigma})=\begin{cases}
\Z_{-}&\textup{if $*=2m-1$}\\
\overline{\ev{\Z/2}}&\textup{if $-1\le *\le 2m-3$ and $*$ is odd}\\
\ev{\Z/2}&\textup{if $0\le *\le 2m-2$ and $*$ is even}
\end{cases}
\end{equation}
The multiplicative generators are also obtained immediately from the CSS and are as follows:
\begin{itemize}
	\item $\zero|u_{\sigma}^{-1}\bar u_{\lambda}^m|\bar u_{\sigma}^{-1}\bar{\bar u}_{\lambda}^m$ generates $\Z_{-}$.
	\item $\zero|u_{\sigma}^{-1}\bar a_{\lambda}^i\bar u_{\lambda}^{m-i}|\zero$ generates $\overline{\ev{\Z/2}}$ %(products in the SS) 
	for $0<i\le m$
	\item $(2a_{\lambda}^iu_{\lambda}^{m-i})/a_{\sigma}|\zero|\zero$ generates $\ev{\Z/2}$ for $0<i\le m$ %(using the $m=1$ case and the fact that $u_{\lambda}$ is an isomorphism since the $2a_{\lambda}u_{\lambda}^k$ generate).
\end{itemize}

%The Kunneth for $H_*(S^{m\lambda-2\sigma})=H_*(S^{m\lambda-\sigma}\wedge S^{\sigma})$ has extension problems, so we instead use the
\subsection{ The case $n=2$}The HSS for $S^{m\lambda-2\sigma}$ similarly collapses with no extensions to give:
\begin{equation}
H_*(S^{m\lambda-2\sigma})=\begin{cases}
\Z&\textup{if $*=2m-2$}\\
\ev{\Z/4}&\textup{if $0\le *<2m-2$ and $*$ is even}\\
Q&\textup{if $*=-2$}
\end{cases}
\end{equation}

%The Kunneth extension problems were
%\begin{gather}
%0\to Q\to \ev{\Z/4}\to \overline{\ev{\Z/2}}\to 0\\
%0\to p^*L\to \Z\to \ev{\Z/2}\to 0
%\end{gather}
The multiplicative generators:
%We note that for $m=0$ the same computation gives $p^*L$ in degree $-2$, as expected. So everything from now on we must require $m\ge 1$.
\begin{itemize}
	\item $u_{\lambda}^m/u_{2\sigma}|u_{\sigma}^{-2}\bar u_{\lambda}^m|\bar u_{\sigma}^{-2}\bar{\bar{u}}_{\lambda}^m$ generates $\Z$. %($2u_{2\sigma}^{-1}u_{\lambda}^m=\Tr(u_{\sigma}^{-2}\bar u_{\lambda}^m)$ but we know this is twice the generator, and we have no torsion!). 
	%The equation
%	\begin{equation}
%	a_{\sigma}u_{2\sigma}^{-1}u_{\lambda}^m=a_{\sigma}^{-1}(2a_{\lambda})u_{\lambda}^{m-1}
%	\end{equation}
%	holds (in $\Z/2$) because multiplying with $a_{\sigma}$ is the injection $H_{2m-2}S^{m\lambda-\sigma}=\ev{\Z/2}\to H_{2m-2}S^{m\lambda}=\ev{\Z/4}$ and then multiplying by $u_{2\sigma}$ is an isomorphism $H_{2m-2}S^{m\lambda}\to H_{2m+2}S^{2\sigma+m\lambda}$ (so in total an injection), and then we get the Gold Relation. 
	%We won't write $a_{\sigma}^{-2}(2a_{\lambda})$ because $a_{\sigma}$ is not an isomorphism $\Z\to \ev{\Z/2}$, but you can think of it as the mod $2$ reduction of $u_{2\sigma}^{-1}u_{\lambda}^m$.
	\item $(a_{\lambda}^iu_{\lambda}^{m-i})/u_{2\sigma}|u_{\sigma}^{-2}\bar a_{\lambda}^i\bar u_{\lambda}^{m-i}|\zero$ generates $\ev{\Z/4}$ for $0<i<m$. %(it's half of $\Tr(u_{\sigma}^{-2}\bar a_{\lambda}^i\bar u_{\lambda}^{m-i})$ and there is no torsion!). 
	%The equation
%	\begin{equation}
%	a_{\sigma}u_{2\sigma}^{-1}u_{\lambda}^{m-i}a_{\lambda}^i=a_{\sigma}^{-1}(2a_{\lambda}^{i+1})u_{\lambda}^{m-i-1}
%	\end{equation}
%	holds (in $\Z/2$) for the same reason. %with the same caveat.
	\item $(2a_{\lambda}^m)/u_{2\sigma}|u_{\sigma}^{-2}\bar a_{\lambda}^m|\zero$ generates $Q$.
\end{itemize}

\subsection{ The case $n=3$}
The KSS for $S^{m\lambda-2\sigma}\wedge S^{-\sigma}$ collapses with no extensions to give:
\begin{equation}
H_*(S^{m\lambda-3\sigma})=\begin{cases}
\Z_{-}&\textup{if $*=2m-3$}\\
\ev{\Z/2}&\textup{if $-2\le *<2m-3$ and $*$ is even}\\
\overline{\ev{\Z/2}}&\textup{if $-1\le *<2m-3$ and $*$ is odd}\\
Q&\textup{if $*=-3$}
\end{cases}
\end{equation}

The multiplicative generators:

\begin{itemize}
	\item $\zero|u_{\sigma}^{-3}\bar u_{\lambda}^m|\bar u_{\sigma}^{-3}\bar{\bar u}_{\lambda}^m$ generates $\Z_{-}$.
	\item $(2a_{\lambda}^iu_{\lambda}^{m-i})/(a_{\sigma}u_{2\sigma})|\zero|\zero$ generates $\ev{\Z/2}$ for $0<i\le m$. 
	%by  as follows: The map $u_{2\sigma}:S^{m\lambda-2\sigma}\to S^{m\lambda}$ induces on $H_{-2}\to H_0$ the map $Q\hookrightarrow \Z/4$ which is an isomorphism in $\ev{\Z/2}=\Tor^1(Q,\Z_{-})\to \Tor^1(\ev{\Z/4},\Z_{-})=\ev{\Z/2}$  and thus by the comparison of Kunneth SS's
	%\begin{equation}
	%\Tor^{p,q}(H_*S^{m\lambda-2\sigma},H_*S^{-\sigma})\to 	\Tor^{p,q+2}(H_*S^{m\lambda},H_*S^{-\sigma})
	%\end{equation}
%	which both collapse with no extensions, $u_{2\sigma}$ is an iso $H_{-2}(S^{m\lambda-3\sigma})\to H_0(S^{m\lambda-\sigma})$. For $0\le *\le 2m-4$ even, $u_{2\sigma}$ is an iso $H_*S^{m\lambda-2\sigma}=\ev{\Z/4}\to \ev{\Z/4}=H_*S^{m\lambda}$ and thus by the same argument is an iso for $S^{m\lambda-3\sigma}\to S^{m\lambda-\sigma}$. In those degrees i.e. for $i<m$, 
%	\begin{equation}
%	u_{2\sigma}^{-1}(a_{\sigma}^{-1}(2a_{\lambda}^i)u_{\lambda}^{m-i})=a_{\sigma}^{-1}(2a_{\lambda}^i)(u_{2\sigma}^{-1}u_{\lambda}^{m-i})
%	\end{equation}	
	%	This could also be  since even though $a_{\sigma}^2$ is $\ev{\Z/2}\to \ev{\Z/4}\to \ev{\Z/2}$ which is trivial, we have an extension problem at the last stage $Q\to \ev{\Z/4}\to \ev{\Z/2}$ and we are not seeing $Q$.
	\item $\zero|u_{\sigma}^{-3}\bar a_{\lambda}^i\bar u_{\lambda}^{m-i}|\zero$ generates $\overline{\ev{\Z/2}}$ for $0<i<m$
	\item  $w_3a_{\lambda}^m|u_{\sigma}^{-3}\bar a_{\lambda}^m|\zero$ generates $Q$. %We don't have $w_3,u_{\lambda}$ together since $w_3u_{\lambda}=0$ as $\overline{\ev{\Z/2}}$ is $0$ in top degree.
\end{itemize}

The generator of $\ev{\Z/2}$ is not immediately obtained from the three spectral sequences. Instead, our argument uses the result of the next subsection for $H_*(S^{m\lambda-4\sigma})$, the computation of which does not use the expression for the $\ev{\Z/2}$ generator so our reasoning is not circular. Given the computation of $H_*(S^{m\lambda-4\sigma})$ we have on top level:
\begin{center}\begin{tikzcd}
	H_{2m-4}(S^{m\lambda-4\sigma})\ar[d,equals]\ar[r,"a_{\sigma}"]&H_{2m-4}(S^{m\lambda-3\sigma})\ar[d,equals]\ar[r,"a_{\sigma}"]&H_{2m-4}(S^{m\lambda-2\sigma})\ar[d,equals]\\
	\Z u_{\lambda}^m/u_{2\sigma}^{2}\ar[r,"a_{\sigma}"]&\Z/2\ar[r,"a_{\sigma}"]&\Z/4(a_{\lambda}u_{\lambda}^{m-1})/u_{2\sigma}
	\end{tikzcd}\end{center}
The composite $a_{\sigma}^2$ sends the generator in $\Z$ to twice the generator in $\Z/4$ by the Gold Relation, hence the second map $a_{\sigma}$ must be the canonical inclusion $\Z/2\to \Z/4$; we conclude that the middle generator is $(2a_{\lambda}u_{\lambda}^{m-1})/(a_{\sigma}u_{2\sigma})$. 

Similarly, we can prove that the other $\ev{\Z/2}$'s in $H_*(S^{m\lambda-3\sigma})$ are generated by $(2a_{\lambda}^iu_{\lambda}^{m-i})/(a_{\sigma}u_{2\sigma})$ (but now in the argument above use $\Z/4$ in the place of $\Z$). 

%We also note that $a_{\sigma}^2$ times these generators is $0$ (since multiplying them by $a_{\sigma}^2u_{2\sigma}$ returns $0$ in light of  $2a_{\sigma}=0$). So $(2a_{\lambda})/a_{\sigma}^{3}$ does \emph{not} exist.\medbreak

\subsection{ The case $n=4$}
For $S^{m\lambda-4\sigma}$ the HSS has only one possibly nontrivial differential $\ev{\Z/4}\to \ev{\Z/2}$ and whether it's 0 or not determines whether $H_{-3}$ is $\ev{\Z/2}$ or $0$. But the KSS collapses and gives $H_{-3}=\ev{\Z/2}$, so the aforementioned differential has to be trivial. So now the HSS  gives
\begin{equation}
H_*(S^{m\lambda-4\sigma})=\begin{cases}
\Z&\textup{if $*=2m-4$}\\
\ev{\Z/4}&\textup{if $-2\le *<2m-4$ and $*$ is even}\\
\ev{\Z/2}&\textup{if $*=-3$}\\
Q&\textup{if $*=-4$}
\end{cases}
\end{equation}
\begin{itemize}
	\item $u_{\lambda}^m/u_{2\sigma}^{2}|u_{\sigma}^{-4}\bar u_{\lambda}^m|\bar u_{\sigma}^{-4}\bar{\bar u}_{\lambda}^m$ generates $\Z$.
	\item $(u_{\lambda}^{m-i}a_{\lambda}^i)/u_{2\sigma}^2|\bar u_{\sigma}^{-4}\bar u_{\lambda}^{m-i}\bar a_{\lambda}^i|\zero$ generates $\ev{\Z/4}$ for $0<i<m$.
	\item $(w_3a_{\lambda}^m)/a_{\sigma}|\zero|\zero$ generates $\ev{\Z/2}$.
	\item $(2a_{\lambda}^m)/u_{2\sigma}^2|u_{\sigma}^{-4}\bar a_{\lambda}^m|\zero$ generates $Q$.
\end{itemize}
%We don't have any extra $\ev{\Z/2}$'s because $w_3u_{\lambda}=0$ by the $H_*(S^{m\lambda-3\sigma})$ computation.\bigbreak

\subsection{ The case $n=5$}
For $S^{m\lambda-5\sigma}$ the HSS has only one possibly nontrivial differential (for $m\ge 2$) and comparison with the KSS shows that differential vanishes. In degree $-3$ for $m> 1$ we have an extension of $\ev{\Z/2}$ and $\overline{\ev{\Z/2}}$ in all three spectral sequences, and the answer can be either $Q$ or $\ev{\Z/2}\oplus \overline{\ev{\Z/2}}$. To see which one it is, we use the multiplicative generators:  $\ev{\Z/2}$ is generated on top level by $(w_3a_{\lambda}^m)/a_{\sigma}^2$ while  $\overline{\ev{\Z/2}}$ is generated on middle level by $u_{\sigma}^{-5}\bar a_{\lambda}^{m-1}\bar u_{\lambda}$ and the question is whether we have the equality:
\begin{equation}
\Tr_2^4(u_{\sigma}^{-5}\bar a_{\lambda}^{m-1}\bar u_{\lambda})\stackrel{?}{=}(w_3a_{\lambda}^m)/a_{\sigma}^2
\end{equation} 
The left hand side is computed by Frobenius to be $w_5u_{\lambda}a_{\lambda}^{m-1}$ so we ask if
\begin{equation}
w_5a_{\lambda}^{m-1}u_{\lambda}\stackrel{?}{=}(w_3a_{\lambda}^m)/a_{\sigma}^2
\end{equation} 
Multiplication by $a_{\sigma}^2$ is an isomorphism hence we equivalently want to check
\begin{equation}
a_{\sigma}^2w_5a_{\lambda}^{m-1}u_{\lambda}\stackrel{?}{=}w_3a_{\lambda}^m
\end{equation}
But by the Gold Relation the left hand side is $0$ as $2w_5=0$, while in the right hand side we have the generator of $\Z/2$. This means that the extension has to be trivial i.e. $\ev{\Z/2}\oplus \overline{\ev{\Z/2}}$.

%Interestingly, $\Tor(Q,\Z_{-})=\Tor(\ev{\Z/2}\oplus \overline{\Z/2})$ so the next Kunneths won't be able to see the difference.

For $m=1$ we have an extension of $\ev{\Z/2}$ and $\Z_{-}$ that resolves to $\ev{\Z/2}\oplus \Z_{-}$ for the same reason. The answer is:

\begin{equation}
H_*(S^{m\lambda-5\sigma})=\begin{cases}
\Z_{-}&\textup{if $*=2m-5$ and $m\ge 2$}\\
\overline{\ev{\Z/2}}&\textup{if $-1\le *<2m-5$ and $*$ is odd}\\
\ev{\Z/2}&\textup{if $-4\le *<2m-5$ and $*$ is even}\\
\ev{\Z/2}\oplus \overline{\ev{\Z/2}}&\textup{if $*=-3$ and $m\ge 2$}\\
\ev{\Z/2}\oplus\Z_{-} &\textup{if $*=-3$ and $m=1$}\\
Q&\textup{if $*=-5$}
\end{cases}
\end{equation}
\begin{itemize}
	\item $\zero|u_{\sigma}^{-5}\bar u_{\lambda}^m|\bar u_{\sigma}^{-5}\bar{\bar u}_{\lambda}^m$ generates all instances of $\Z_{-}$.
	\item $\zero|u_{\sigma}^{-5}\bar u_{\lambda}^{m-i}\bar a_{\lambda}^i|\zero$ generates all instances of $\overline{\ev{\Z/2}}$ for $0<i\le m-1$.
	\item $(2a_{\lambda}^iu_{\lambda}^{m-i})/(a_{\sigma}u_{2\sigma}^{2})|\zero|\zero$ generates $\ev{\Z/2}$ for $0<i\le m$
	\item $(w_3a_{\lambda}^m)/a_{\sigma}^2|\zero|\zero$ generates the $\ev{\Z/2}$ that appear as summands at $-3$.
	\item  $w_5a_{\lambda}^m|u_{\sigma}^{-5}\bar a_{\lambda}^m|\zero$ generates $Q$ at $-5$.
\end{itemize} 
The multiplicative generator $(2a_{\lambda}^iu_{\lambda}^{m-i})/(a_{\sigma}u_{2\sigma}^{2})$ is proven just like we did with $H_*(S^{m\lambda-3\sigma})$ (i.e. we use the result for $H_*(S^{m\lambda-6\sigma})$ and multiply by $a_{\sigma}^2$).
%Note that this implies $w_5u_{\lambda}=(w_3/u_{2\sigma})u_{\lambda}=0$.
\medbreak

\subsection{ The general case} We proceed by induction on $n$.

 The computation of $H_*(S^{m\lambda-n\sigma})$ for even $n\ge 6$ is exactly like that for $n=4$: The HSS has unknown differentials that are trivial by comparison with the KSS and there are no extension problems.
 
 The computation of $H_*(S^{m\lambda-n\sigma})$ for odd $n\ge 5$ is exactly like that for $n=5$: The differentials in the HSS vanish by comparison with the KSS and we have extension problems for $\overline{\ev{\Z/2}},\ev{\Z/2}$ and $\Z_{-},\ev{\Z/2}$ that both resolve to trivial extensions for the same reason as the $n=5$ case. 
 
 The results are displayed in subsections \ref{Lambda minus sigma even} and \ref{Lambda minus sigma odd}.

%\stepcounter{proofpart}
\section{Preparation for \texorpdfstring{$H_*(S^{n\sigma-m\lambda})$}{the case n*sigma-m*lambda}}

In this section and the next we will compute $H_*(S^{n\sigma-m\lambda})$ for $n,m\ge 1$. As before, there are three spectral sequences of Mackey functors we will use, and this section is devoted to determining their $E_2$ terms sans differentials. We omit the multiplicative presentation for certain generators as in subsection \ref{ProofsE2page1}.

\subsection{The Homological Spectral Sequence} The HSS is
\begin{equation}
E^2_{p,q}=H_p(S^{n\sigma};H_{q}S^{-m\lambda})\implies H_{p+q}(S^{n\sigma-m\lambda})
\end{equation}
For $n\ge 1$,
\begin{gather}
H_*(S^{n\sigma};H_{-2m}(S^{-m\lambda}))=H_*(S^{n\sigma};L)=\begin{cases}
L^{\sharp}&\textup{if $*=n$: even}\\
\Z_{-}^{\flat}&\textup{if $*=n$: odd}\\
\ev{\Z/2}&\textup{if $2\le *<n$ and $*$ is even} \end{cases}
\end{gather}
The multiplicative generators:
\begin{itemize}
	\item $(2u_{2\sigma}^{n/2})/u_{\lambda}^{m}|(2u_{\sigma}^n)/\bar u_{\lambda}^{m}|\bar u_{\sigma}^n\bar{\bar{u}}_{\lambda}^{-m}$ generates $L^{\sharp}$.
	\item $\zero|(2u_{\sigma}^n)/\bar u_{\lambda}^{m}|\bar u_{\sigma}^n\bar{\bar{u}}_{\lambda}^{-m}$ generates $\Z_{-}^{\flat}$
\end{itemize}
\phantom{}\medbreak

For $1\le i\le m-1$,
\begin{gather}
H_*(S^{n\sigma};H_{-2i-1}(S^{-m\lambda}))=H_*(S^{n\sigma};\ev{\Z/4})=\begin{cases}
Q^{\sharp}&\textup{if $*=n\neq 1$}\\
\ev{\Z/2}&\textup{if $0\le *<n$ and $*\neq 1$}\\
\overline{\ev{\Z/2}}&\textup{if $*=n=1$}\end{cases}
\end{gather}
The multiplicative generators:
\begin{itemize}
	\item $(u_{2\sigma}^{n/2}s)/(a_{\lambda}^{m-i}u_{\lambda}^{i-2})|(u_{\sigma}^n\bar s)/(\bar a_{\lambda}^{m-i}\bar u_{\lambda}^{i-2})|\zero$ generates $Q^{\sharp}$ for $*=n$ even.
	\item  The middle level generator of $Q^{\sharp}$ for $*=n$ odd is $(u_{\sigma}^n\bar s)/(\bar a_{\lambda}^{m-i}\bar u_{\lambda}^{i-2})$.
	\item  $\zero|(u_{\sigma}\bar s)/(\bar a_{\lambda}^{m-i}\bar u_{\lambda}^{i-2})|\zero$ generates $\overline{\ev{\Z/2}}$
	\item  $(a_{\sigma}^{n-*}u_{2\sigma}^{*/2}s)/(a_{\lambda}^{m-i}u_{\lambda}^{i-2})|\zero|\zero$ generates $\ev{\Z/2}$ for $0\le *<n$ even.
\end{itemize}

\phantom{}\medbreak
\subsection{The Cohomological  Spectral Sequence} The CSS is
\begin{equation}
E_2^{p,q}=H_{-p}(S^{-m\lambda};H_{-q}S^{n\sigma})\implies H_{-p-q}(S^{n\sigma-m\lambda})
\end{equation}
and for odd $n$,
\begin{gather}
H_{-*}(S^{m\lambda};H^{n}(S^{n\sigma}))=H_{-*}(S^{m\lambda};\Z_{-})=\begin{cases}
L_{-}&\textup{if $*=2m$}\\
\overline{\ev{\Z/2}}&\textup{if $3\le *<2m$ and $*$ is odd}\\
\ev{\Z/2}&\textup{if $2\le *<2m$ and $*$ is even}
\end{cases}
\end{gather}
The multiplicative generators:
\begin{itemize}
	\item $\Tr_2^4((2u_{\sigma}^n)/\bar u_{\lambda}^{m})|(2u_{\sigma}^n)/\bar u_{\lambda}^{m}|\bar u_{\sigma}^n\bar{\bar{u}}_{\lambda}^{-m}$ generates $L_{-}$
	\item $\zero|(u_{\sigma}^n\bar s)/(\bar a_{\lambda}^{*-2}\bar u_{\lambda}^{m-*})|\zero$ generates $\overline{\ev{\Z/2}}$ 
\end{itemize}
\phantom{}\medbreak

Finally, $H_{-*}(S^{-m\lambda};\ev{\Z/2})=\ev{\Z/2}$ concentrated in degree $0$. So the homology  $H_{-*}(S^{-m\lambda};H_iS^{n\sigma})$ is generated by
\begin{equation}(a_{\sigma}^{n-2i}u_{2\sigma}^i)/a_{\lambda}^m\end{equation}
for $0\le i< n/2$.\medbreak
%(The map $a_{\lambda}:C^*(S^{m\lambda};H_iS^{n\sigma})\to C^*(S^0;H_iS^{n\sigma})$ is projection).\medbreak

\subsection{The Kunneth Spectral Sequence} 

The KSS is
\begin{equation}
E^2_{p,q}=\Tor^{p,q+2}_{\underline{\Z}}(H_*S^{n\sigma-(m-1)\lambda},L)\implies H_{p+q}(S^{n\sigma-m\lambda})
\end{equation}
This is the Kunneth spectral sequence for $S^{n\sigma-(m-1)\lambda}\wedge S^{-\lambda}$ and we have used that $H_*S^{-\lambda}=L$ concentrated in degree $*=-2$. A free resolution of $L$ is
\begin{gather}
%0\to \Z\xrightarrow{1\mapsto \sum g^ix} \Z_{C_4}\xrightarrow{x\mapsto x-gx} \Z_{C_4}\xrightarrow{x\mapsto 1} \Z \to \ev{\Z/4}\to 0\\
%0\to \Z\xrightarrow{x-gx} \Z_{C_2}\xrightarrow{x\mapsto 1} \Z \to \ev{\Z/2}\to 0\\
0\to \Z\xrightarrow{\sum g^ix} \Z_{C_4}\xrightarrow{x\mapsto x-gx} \Z_{C_4}\xrightarrow{x\mapsto 1} L\to 0
\end{gather}
where $x\in \Z[C_4]$ corresponds to a generator of $C_4$ (recall that $\Z_{C_4}$ is the fixed point Mackey functor on $\Z[C_4]$). 

We list only the nonzero $\Tor$'s below, dropping the subscript $\underline{\Z}$ from the notation:
\begin{itemize}
	\setlength\itemsep{0.5em}
%	\item $\Tor(\ev{\Z/4},\Z_{-})=\overline{\ev{\Z/2}},\ev{\Z/2},0,...$ (checks out)
%	\item $\Tor(\ev{\Z/4},\ev{\Z/4})=\ev{\Z/4},0,0,\ev{\Z/4},0,...$
%	\item $\Tor(\ev{\Z/4},\ev{\Z/2})=\ev{\Z/2},0,0,\ev{\Z/2},0,...$ (checks out)
%%	\item $\Tor(\ev{\Z/4},\overline{\ev{\Z/2}})=\overline{\ev{\Z/2}},0,0,\overline{\ev{\Z/2}},0,...$ (checks out)
%	\item $\Tor(\ev{\Z/2},\Z_{-})=0,\ev{\Z/2},0,...$
%	\item $\Tor(\ev{\Z/2},\ev{\Z/2})=\ev{\Z/2},0,\ev{\Z/2},0,...$
	\item $\Tor^0(L,\Z_{-})=L_{-}$ %(checks out)
	\item $\Tor^2(L,\ev{\Z/4})=\ev{\Z/4}$
	\item $\Tor^2(L,\ev{\Z/2})=\ev{\Z/2}$ %(checks out)
	\item $\Tor^0(L,L_{-})=L_{-}$ and $\Tor^1(L,L_{-})=\overline{\ev{\Z/2}}$ and $\Tor^2(L,L_{-})=\ev{\Z/2}$
	\item $\Tor^2(L,\overline{\ev{\Z/2}})=\overline{\ev{\Z/2}}$
	\item $\Tor^0(L,L)=L$ and $\Tor^1(L,L)=\ev{\Z/4}$ %(checks out)
	\item $\Tor^0(L,L^{\sharp})=L$ and $\Tor^1(L,L^{\sharp})=Q^{\sharp}$.
	\item $\Tor^0(L,\Z_{-}^{\flat})=L_{-}$ and $\Tor^1(L,\Z_{-}^{\flat})=\overline{\ev{\Z/2}}$.
	\item $\Tor^2(L,Q^{\sharp})=Q^{\sharp}$.
\end{itemize}
 The existence of $\Tor^2$ terms means that the KSS can now have potentially non vanishing differentials and more complicated extension problems. As a result the computations in the next section are slightly more involved compared to those in section \ref{Proofssigma}.

Still, this is slightly better than the worst case scenario that is nonvanishing $\Tor^3$: In general, for finite cyclic $G$ the abelian category of $\underline{\Z}$-Mackey functors has projective dimension $3$ (\cite{BSW17}).

%\stepcounter{proofpart}
\section{Proofs for \texorpdfstring{$H_*(S^{n\sigma-m\lambda})$}{The case n*sigma-m*lambda}}\label{Proofslambda}
The computation of  $H_*(S^{n\sigma-m\lambda})$ depends heavily on the parity of $n$ and we distinguish three cases: even $n$, $n=1$ and odd $n\ge 3$.

In the even $n$ case, we compute $H_*(S^{n\sigma-m\lambda})$ for $m=1,2,3$ separately before we can perform induction. The case $n=1$ is staightforward enough to do for all $m$ at once, while for odd $n\ge 3$ we again compute the special cases $m=1,2,3$ and then induct on $m$. 

We make use of all three spectral sequences HSS, CSS and KSS and play them off against each other.

\subsection{The case of even $n$ and $m=1$}
The HSS for $S^{n\sigma-\lambda}$ collapses with no extensions (for dimensional reasons) to give
\begin{equation}
H_*(S^{n\sigma-\lambda})=
\begin{cases}L^{\sharp}&\textup{if  $*=n-2$}\\
\ev{\Z/2}&\textup{if $0\le *\le n-4$ and $*$ is even}
\end{cases}
\end{equation}
The multiplicative generators for $L^{\sharp}$ and $\ev{\Z/2}$ are obtained immediately from the HSS and CSS respectively and they are:
\begin{itemize}
	\item  $(2u_{2\sigma}^{n/2})/u_{\lambda}|(2u_{\sigma}^n)/\bar u_{\lambda}|\bar u_{\sigma}^n\bar{\bar{u}}_{\lambda}^{-1}$ generates $L^{\sharp}$.
	\item  $(a_{\sigma}^{n-2i}u_{2\sigma}^i)/a_{\lambda}|\zero|\zero$ generates $\ev{\Z/2}$ for $0\le i<n/2-1$%(or the Kunneth argument: $a_{\lambda}$ is an iso for $S^{n\sigma}\to S^{n\sigma+\lambda}$ in $H_*\to H_{*+2}$ except for $*=n$ where it's $\ev{\Z/2}\hookrightarrow \ev{\Z/4}$, so it induces an iso in the Kunneth's $S^{n\sigma}\wedge S^{-\lambda}\to S^{n\sigma+\lambda}\wedge S^{-\lambda}$ and there are no extensions in the relevant part).
\end{itemize}
%The Gold Relation implies that the $\ev{\Z/2}$ generators multiplied with $u_{\lambda}$ give $0$ ($2a_{\sigma}$).\medbreak

Note that by the Gold relation, the mod $2$ reduction of $(2u_{2\sigma}^i)/u_{\lambda}$ is $(a_{\sigma}^2u_{2\sigma}^{i-1})/a_{\lambda}$. In particular, $(2u_{2\sigma})/u_{\lambda}$ and its mod $2$ reduction $a_{\sigma}^2/a_{\lambda}$ exist. However as we will see, $(2u_{2\sigma})/u_{\lambda}^2$ and $a_{\sigma}/a_{\lambda}$ do \emph{not} exist.

\subsection{The case of even $n$ and $m=2$}
The KSS for $S^{n\sigma-2\lambda}$ collapses %(alternatively, compare the two AHSS's) 
giving the answer with the exception of degree $n-4$ for $n\ge 4$. There is an extension problem of $L,\ev{\Z/2}$ in all three spectral sequences and there are only two possibilities for that extension:
\begin{gather}
0\to L\to L\oplus \ev{\Z/2}\to \ev{\Z/2}\to 0\\
0\to L\to L^{\sharp}\to \ev{\Z/2}\to 0
\end{gather}
By the CSS, the extension is $L^{\sharp}$ iff $(2u_{2\sigma}^{n/2})/u_{\lambda}^{2}$ exists and its mod $2$ reduction is $(a_{\sigma}^4u_{2\sigma}^{n/2-2})/a_{\lambda}^2$. If this were true, then  we would have the following commutative diagram:

\begin{center}
	\begin{tikzcd}
&& \Z (u_{2\sigma}^{n/2})\ar[rd,"a_{\lambda}^2"]&\\
&\Z \left(\frac{ 2u_{2\sigma}^{n/2}}{u_{\lambda}^2}\right)\ar[dr,"a_{\lambda}^2"]\ar[ur,"u_{\lambda}^2"]%\ar[dl,"\text{mod 2}" left]
&& \Z/4 (u_{2\sigma}^{n/2}a_{\lambda}^2)\\
%\Z/2  \left(\frac{a_{\sigma}^4u_{2\sigma}^{n/2-2}}{a_{\lambda}^2}\right)\ar[rr,"a_{\lambda}^2"]
&&\Z/2  \left(a_{\sigma}^4u_{2\sigma}^{n/2-2}\right)\ar[ru,"u_{\lambda}^2"]&
	\end{tikzcd}
\end{center}
Note that in the lower part of the diagram, since $a_{\sigma}^3u_{\lambda}=0$ by the Gold Relation, the $u_{\lambda}^2$ map is trivial. So we have
\begin{center}
	\begin{tikzcd}
	& \Z\ar[rd,"\text{mod 4}"]&\\
	\Z\ar[dr,"\text{mod 2}"]\ar[ur,"2"]&& \Z/4\\
	&\Z/2 \ar[ru,"0"]&
	\end{tikzcd}
\end{center}
which clearly doesn't commute.

Therefore the extension has to be $L\oplus \ev{\Z/2}$, which means that $(2u_{2\sigma}^{n/2})/u_{\lambda}^2$ does not exist. As we remarked above, this only happens for $n\ge 4$; for $n=2$ we only have an $L$ so there is no extension and the elements $(2u_{2\sigma})/u_{\lambda}^{2}$ and $a_{\sigma}^2/a_{\lambda}^{2}$ do \emph{not} exist. 

\medbreak
In conclusion we have
\begin{equation}
H_*(S^{n\sigma-2\lambda})=
\begin{cases}Q^{\sharp}&\textup{if  $*=n-3$}\\
L\oplus \ev{\Z/2}&\textup{if  $*=n-4$ and $n\ge 4$}\\
L&\textup{if  $*=-2$ and $n=2$}\\
\ev{\Z/2}&\textup{if $0\le *<n-4$ and $*$ is even}\\
\end{cases}
\end{equation}

\begin{itemize}
	\item  $u_{2\sigma}^{n/2}s|u_{\sigma}^n\bar s|\zero$ generates $Q^{\sharp}$ %(by the cohomological)% or the Kunneth for $S^{n\sigma}\wedge S^{-2\lambda}$).
	\item  $(4u_{2\sigma}^{n/2})/u_{\lambda}^{2}|(2u_{\sigma}^n)/\bar u_{\lambda}^{2}|\bar u_{\sigma}^{n}\bar{\bar{u}}_{\lambda}^{-2}$ generates all instances of $L$ (both as a summand and nonsummand).
	\item  	 $(a_{\sigma}^{n-2i}u_{2\sigma}^i)/a_{\lambda}^{2}|\zero|\zero$ generates all instances of $\ev{\Z/2}$ for $0\le i<n/2-1$ %(by the cohomological).
	
\end{itemize}

\subsection{The case of even $n$ and $m=3$}

The KSS for $S^{n\sigma-3\lambda}$ collapses with an usual extension problem of $L$ and $\ev{\Z/2}$ at $n-6$ and $n\ge 6$. The answer is $L^{\sharp}$ iff $(2u_{2\sigma}^{n/2})/u_{\lambda}^{3}$ exists, but if it did then $(2u_{2\sigma}^{n/2})/u_{\lambda}^2$ would also exist, contradicting the computation of $S^{n\sigma-2\lambda}$ in the preceding subsection. We conclude:
\begin{equation}
H_*(S^{n\sigma-3\lambda})=
\begin{cases}Q^{\sharp}&\textup{if  $*=n-3$}\\
\ev{\Z/4}&\textup{if  $*=n-5$}\\
L\oplus \ev{\Z/2}&\textup{if $*=n-6$ and $n\ge 6$}\\
L&\textup{if $*=n-6$ and $n=2,4$}\\
\ev{\Z/2}&\textup{if $0\le *\le n-4$ and $*$ is even and $*\neq n-6$}
\end{cases}
\end{equation}

\begin{itemize}
	\item $(u_{2\sigma}^{n/2}s)/a_{\lambda}|(u_{\sigma}^n\bar s)/\bar a_{\lambda}|\zero$ generates $Q^{\sharp}$ 
	\item $u_{2\sigma}^{n/2}(s/u_{\lambda})|(u_{\sigma}\bar s)/\bar u_{\lambda}|\zero$ generates $\ev{\Z/4}$
	\item $(a_{\sigma}^{2i}u_{2\sigma}^{n/2-i})/a_{\lambda}^{3}|\zero|\zero$ generates all instances of $\ev{\Z/2}$ for $2\le i\le n/2$.
	\item $(4u_{2\sigma}^{n/2})/u_{\lambda}^3|(2u_{\sigma}^n)/\bar u_{\lambda}^3|\bar u_{\sigma}^{n}\bar{\bar u}_{\lambda}^{-3}$ generates all instances of $L$.
\end{itemize}

\subsection{The general case of even $n$}

We proceed by induction, with the case of $H_*S^{n\sigma-m\lambda}$ for $m\ge 4$ being treated exactly the same as for $m=3$. The answer is given in subsection \ref{Sigma minus Lambda even}. 
\iffalse
For even $n$ by induction we get:
\begin{equation}
H_*(S^{n\sigma-m\lambda})=
\begin{cases}Q^{\sharp}&\textup{if  $*=n-3$ and $m\ge 2$}\\
\ev{\Z/4}&\textup{if  $n-2m<*<n-3$ and $*$ is odd}\\
\ev{\Z/2}&\textup{if $0\le *\le n-4$ and $*$ is even and $*\neq n-2m$}\\
L\oplus \ev{\Z/2}&\textup{if $*=n-2m$ and $n-2m\ge 0$ and $m\ge 2$}\\
L&\textup{if $*=n-2m$ and $n-2m<0$ and $m\ge 2$}\\
L^{\sharp}&\textup{if $*=n-2$ and $m=1$}
\end{cases}
\end{equation}

Note: For $m=2$, MATLAB may return the extension of $L,\ev{\Z/2}$ that has transfer $\Tr:\Z\to \Z\oplus \Z/2$ equal to $\Tr(1)=(1,1)$. This Mackey functor is isomorphic to $L\oplus \ev{\Z/2}$ through $\Z/2\oplus \Z\to \Z/2\oplus \Z$ that is $(1,0)\mapsto (1,0)$ and $(0,1)\mapsto (1,1)$.

\begin{itemize}
	\item $u_{2\sigma}^{n/2}a_{\lambda}^{-m+2}s$ generates $Q^{\sharp}$ 
	\item $u_{2\sigma}^{n/2}a_{\lambda}^{-(i-2)}u_{\lambda}^{-(m-i)}s$ generates $\ev{\Z/4}$ for $2\le i<m$
	\item $a_{\lambda}^{-m}a_{\sigma}^{2i}u_{2\sigma}^{n/2-i}$ generates $\ev{\Z/2}$ for $2\le i\le n/2$.
	\item $4u_{\lambda}^{-m}u_{2\sigma}^{n/2}$ generates $L$.
	\item  $2u_{\lambda}^{-1}u_{2\sigma}^{n/2}$ generates $L^{\sharp}$ for $m=1$
\end{itemize}
\fi

\subsection{The case of $n=1$}

By comparing the HSS (collapses with one extension) and CSS (has differentials but no extensions) we get:

\begin{equation}
H_*(S^{\sigma-m\lambda})=
\begin{cases}L_{-}&\textup{if  $*=1-2m$ and $m\ge 2$}\\
\Z_{-}^{\flat}&\textup{if $*=-1$ and $m=1$}\\
\overline{\ev{\Z/2}}&\textup{if $-2m+2\le *\le -2$ and $*$ is even}\\
\ev{\Z/2}&\textup{if $-2m+3\le *\le -3$ and $*$ is odd}
\end{cases}
\end{equation}

\begin{itemize}
	\item $(a_{\sigma}s)/u_{\lambda}^{m-2}|(2u_{\sigma})/\bar u_{\lambda}^m|\bar u_{\sigma}\bar{\bar{u}}_{\lambda}^{-m}$ generates $L_{-}$ for $m\ge 2$ % (by the homological). 
	\item $\zero|(2u_{\sigma})/\bar u_{\lambda}|\bar u_{\sigma}\bar{\bar{u}}_{\lambda}^{-1}$ generates $\Z_{-}^{\flat}$ for $m=1$
	\item $\zero|(u_{\sigma}\bar s)/(\bar a_{\lambda}^{m-i}\bar u_{\lambda}^{i-2})|\zero$ generates $\overline{\ev{\Z/2}}$ for $2\le i\le m$.
	\item $(a_{\sigma}s)/(a_{\lambda}^{m-i}u_{\lambda}^{i-2})|\zero|\zero$ generates $\ev{\Z/2}$ for $2\le i<m$% (by the homological).
\end{itemize}
\medbreak

Therefore,
\begin{equation}
\Tr_2^4((2u_{\sigma})/\bar u_{\lambda}^2)=a_{\sigma}s
\end{equation}

We also note that while $a_{\sigma}/a_{\lambda}$ does not exist, the element $a_{\sigma}^2/a_{\lambda}$ does exist by the $n=2$ computation. %This element can be understood geometrically as the $C_4$-map $S^{\lambda}\to S^{2\sigma}$ given by $z\mapsto z^2$ on the representations.

\subsection{The case of odd $n\ge 3$ and $m=1$}

The HSS for $S^{n\sigma-\lambda}$ collapses with no extensions to give:
\begin{equation}
H_*(S^{n\sigma-\lambda})=
\begin{cases}\Z_{-}^{\flat}&\textup{if $*=n-2$}\\
\ev{\Z/2}&\textup{if $0\le *\le n-3$ and $*$ is even}
\end{cases}
\end{equation}

\begin{itemize}
	\item $(2u_{\sigma}^n)/\bar u_{\lambda}|\bar u_{\sigma}^n\bar{\bar{u}}_{\lambda}^{-1}|\zero$ generates $\Z_{-}^{\flat}$.
	\item  $(a_{\sigma}^{n-2i}u_{2\sigma}^i)/a_{\lambda}|\zero|\zero$ generates $\ev{\Z/2}$ for $0\le i\le (n-3)/2$ %(by the cohomological). %(or the Kunneth argument: $a_{\lambda}$ is an iso for $S^{n\sigma}\to S^{n\sigma+\lambda}$ in $H_*\to H_{*+2}$ except for $*=n$ where it's $\ev{\Z/2}\hookrightarrow \ev{\Z/4}$, so it induces an iso in the Kunneth's $S^{n\sigma}\wedge S^{-\lambda}\to S^{n\sigma+\lambda}\wedge S^{-\lambda}$ and there are no extensions in the relevant part).
\end{itemize}	

\subsection{The case of odd $n\ge 3$ and $m=2$}
For $S^{n\sigma-2\lambda}$ the HSS and CSS comparison reveals:

\begin{equation}
H_*(S^{n\sigma-2\lambda})=
\begin{cases}Q^{\sharp}&\textup{if  $*=n-3$}\\
L_{-}&\textup{if  $*=n-4$}\\
\ev{\Z/2}&\textup{if $0\le *\le n-5$ and $*$ is even}\\
\end{cases}
\end{equation}

\begin{itemize}
	\item $(a_{\sigma}^3u_{2\sigma}^{(n-3)/2})/a_{\lambda}^2|u_{\sigma}^n\bar s|\zero$ generates $Q^{\sharp}$ % (by the cohomological). 
	\item $a_{\sigma}u_{2\sigma}^{(n-1)/2}s|(2u_{\sigma}^n)/\bar u_{\lambda}^2|\bar u_{\sigma}^n\bar{\bar{u}}_{\lambda}^{-2}$ generates $L_{-}$ % (by the homological).
	\item  	 $(a_{\sigma}^{n-2i}u_{2\sigma}^i)/a_{\lambda}^2|\zero|\zero$ generates $\ev{\Z/2}$ for $0\le i\le (n-5)/2$ %(by the cohomological).
\end{itemize}
So $a_{\sigma}^3/a_{\lambda}^2$ exists and
\begin{equation}
\bar s=\Res^4_2(a_{\sigma}^3/a_{\lambda}^2)/u_{\sigma}^3
\end{equation}

\subsection{The case of odd $n\ge 3$ and $m=3$}
For $S^{n\sigma-3\lambda}$ comparison of the KSS and HSS gives the answer with the exception of an extension problem of $\overline{\ev{\Z/2}}$ and $\ev{\Z/2}$ for $n\ge 5$. There are two possible extensions, $Q^{\sharp}$ and $\overline{\ev{\Z/2}}\oplus \ev{\Z/2}$, and to determine which one it is, we use the multiplicative generators: The middle level generator of $\overline{\ev{\Z/2}}$ is $u_{\sigma}^n(\bar s/\bar u_{\lambda})$ and the top level generator of $\ev{\Z/2}$ is $(a_{\sigma}^5u_{2\sigma}^{(n-5)/2})/a_{\lambda}^3$ so it all rests on whether or not
\begin{equation}
\Res^4_2(a_{\sigma}^5/a_{\lambda}^3)\stackrel{?}{=}u_{\sigma}^5(\bar s/\bar u_{\lambda})
\end{equation}
But we already know that $a_{\sigma}^3/a_{\lambda}^2$ and $a_{\sigma}^2/a_{\lambda}$ both exist, the latter generating the top level of a $\ev{\Z/2}$ thus having trivial restriction. Therefore
\begin{equation}
\Res^4_2(a_{\sigma}^5/a_{\lambda}^3)=\Res^4_2(a_{\sigma}^3/a_{\lambda}^2)\Res^4_2(a_{\sigma}^2/a_{\lambda})=0
\end{equation}
can't be the generator $u_{\sigma}^5(\bar s/\bar u_{\lambda})$ and the extension is trivial. We conclude:

\begin{equation}
H_*(S^{n\sigma-3\lambda})=
\begin{cases}Q^{\sharp}&\textup{if  $*=n-3$}\\
\ev{\Z/2}&\textup{if $*=n-4$}\\
\overline{\ev{\Z/2}}\oplus \ev{\Z/2}&\textup{if $*=n-5$ and $n\ge 5$}\\
\overline{\ev{\Z/2}}&\textup{if $*=n-5$ and $n=3$}\\
L_{-}&\textup{if  $*=n-6$}\\
\ev{\Z/2}&\textup{if $0\le *\le n-7$ and $*$ is even}\\
\end{cases}
\end{equation}

\begin{itemize}
	\item $(a_{\sigma}^3u_{2\sigma}^{(n-3)/2})/a_{\lambda}^3|(u_{\sigma}^n\bar s)/\bar a_{\lambda}|\zero$ generates $Q^{\sharp}$. % (by the cohomological).
	\item $(a_{\sigma}u_{2\sigma}^{(n-1)/2}s)/a_{\lambda}|\zero|\zero$ generates $\ev{\Z/2}$ at degree $n-4$. % (by the homological).
	\item $(a_{\sigma}^{n-2i}u_{2\sigma}^i)/a_{\lambda}^3|\zero|\zero$ generates all instances of $\ev{\Z/2}$ for $0\le i\le (n-5)/2$.
	\item $\zero|(u_{\sigma}^n\bar s)/\bar u_{\lambda}|\zero$ generates all instances of $\overline{\ev{\Z/2}}$.
	\item $(a_{\sigma}u_{2\sigma}^{(n-1)/2}s)/u_{\lambda}|(2u_{\sigma}^n)/\bar u_{\lambda}^3|\bar u_{\sigma}^n\bar{\bar{u}}_{\lambda}^{-3}$ generates $L_{-}$.% (by the homological).
\end{itemize}

\subsection{The general case of odd $n\ge 3$}

We proceed by induction, with $H_*S^{n\sigma-m\lambda}$ for $m\ge 4$ being treated exactly like the $m=3$ case. The answer is given in subsection \ref{Sigma minus Lambda odd}. 

\appendix

\section{The relations}\label{Appendix}

In this appendix we prove that the secondary relations (the Frobenius relations combined with the additive structure and the presentation of the generators given in section \ref{Results}) and the four "extra" relations
\begin{gather}
a_{\sigma}^2u_{\lambda}=2u_{2\sigma}a_{\lambda} \label{Gold}\\
\frac{x_{1,1}}{a_{\sigma}^2}=\frac{w_3}{a_{\lambda}} \label{xw}\\
\frac{x_{1,1}}{a_{\lambda}}=\frac{2s}{a_{\sigma}} \label{xs1}\\
\frac{x_{1,1}}{u_{\lambda}}=\frac{a_{\sigma}s}{u_{2\sigma}} \label{xs2}
\end{gather}
can be used to generate all other relations in $\underline \pi_{\bigstar}^{C_4}(H\underline \Z)$. The first extra relation is the Gold relation and the final two follow from the definition of $s$. The second is actually redundant (we only use it as a convenient way to pass between $x_{1,1}$ and $w_3$) and follows from the Gold in this way: First,
\begin{equation}\label{Eq1}\frac{w_3}{a_{\sigma}^2}\cdot \frac{2u_{2\sigma}}{u_{\lambda}}=\frac{x_{1,1}}{a_{\sigma}^2}\end{equation}
To see this, note that multiplication by $a_{\sigma}^2$ is an isomorphism so equivalently:
\begin{equation}
w_3 \frac{2u_{2\sigma}}{u_{\lambda}}=x_{1,1}
\end{equation}
This is proven by appealing to the Frobenius relation:
\begin{equation}
w_3 \frac{2u_{2\sigma}}{u_{\lambda}}=\Tr_2^4(u_{\sigma}^{-3})\frac{2u_{2\sigma}}{u_{\lambda}}=\Tr_2^4\Big(u_{\sigma}^{-3}\frac{2u_{\sigma}^2}{\bar u_{\lambda}}\Big)=\Tr_2^4\Big(u_{\sigma}^{-1}\frac{2}{\bar u_{\lambda}}\Big)=\Tr_1^4(\bar u_{\sigma}^{-1}\bar{\bar u}_{\lambda}^{-1})=x_{1,1}\\
\end{equation}
Next note that
\begin{equation}\label{Eq2}a_{\lambda}\cdot \frac{2u_{2\sigma}}{u_{\lambda}}=a_{\sigma}^2\end{equation}
as multiplying by $u_{\lambda}$ is an isomorphism (the map $\Z/2\to \Z/2$, $a_{\sigma}^2\mapsto a_{\sigma}^2u_{\lambda}=2u_{2\sigma}a_{\lambda}$ is an isomorphism). By \eqref{Eq1} and \eqref{Eq2},
\begin{equation}
w_3=\frac{x_{1,1}}{a_{\sigma}^2}a_{\lambda}
\end{equation}
So the map $a_{\lambda}$ is an isomorphism $\Z/2\to \Z/2$, $x_{1,1}/a_{\sigma}^2\mapsto w_3$, hence we can write
\begin{equation}
\frac{x_{1,1}}{a_{\sigma}^2}=\frac{w_3}{a_{\lambda}}
\end{equation}

To prove that the secondary relations and the four extra relations are enough to generate all others, it is enough to compute the product of any two generators $a\in H_k^{C_4}(S^{n\sigma+m\lambda})$ and $b\in H_{k'}^{C_4}(S^{n'\sigma+m'\lambda})$ as a linear combination of the generators in $H_{k+k'}^{C_4}(S^{(n+n')\sigma+(m+m')\lambda})$, using only the relations above.

What follows is an exhaustive list of all these products that need to be computed and the results of the computations. The proofs are rather brief; consult subsection \ref{Subtle} for the strategy employed. To keep the length of the list reasonable, we have made the following omissions:

\begin{itemize}
	\item[1] We omit products where one factor is a transfer, as these reduce to the $C_2$ case by the Frobenius relation:
	\begin{equation}
\Tr_2^4(x)y=	\Tr_2^4(x\Res^4_2y)
	\end{equation}
	\item[2] We omit products that are trivial for degree reasons.
	\item [3] We omit products where both factors are in $H_kS^{n\sigma+m\lambda}$ for $k,n,m\ge 0$. This part is polynomially generated by the Euler and orientation classes modulo the Gold relation \eqref{Gold}.
	\item [4] We omit products that can immediately be computed through the following fact: If $x/y, z/w$ and $(xz)/(yw)$ all generate the homology groups they live in, then
	\begin{equation}
	\frac{x}{y}\cdot \frac{z}{w}=\frac{xz}{yw}
	\end{equation}
	Note that this applies only when we have cyclic homology in the degrees of  $x/y, z/w$ and $(xz)/(yw)$.
\end{itemize}

With all that said, we are ready to present the list (we only label the relations/equations that we reference later in the proofs of other relations):
\begin{itemize}	\setlength\itemsep{0.8em}
	\item The following relations compute the product of $a_{\sigma}$ with the other generators:
	\begin{gather}
	a_{\sigma}\cdot \dfrac{s}{u_{2\sigma}^iu_{\lambda}^j}=\frac{x_{1,1}}{u_{2\sigma}^{i-1}u_{\lambda}^{j+1}}\\
	a_{\sigma}\cdot \dfrac{s}{u_{2\sigma}^ia_{\lambda}^ju_{\lambda}^k}=\frac{2s}{a_{\sigma}u_{2\sigma}^{i-1}a_{\lambda}^{j-1}u_{\lambda}^{k+1}}\\
	a_{\sigma}\cdot \frac{x_{1,1}}{u_{2\sigma}^ia_{\lambda}^j}=\frac{2s}{u_{2\sigma}^ia_{\lambda}^{j-1}}\\
	a_{\sigma}\cdot 	\dfrac{u_{\lambda}^i}{u_{2\sigma}^j}=\frac{2a_{\lambda}u_{\lambda}^i}{a_{\sigma}u_{2\sigma}^{j-1}}\\
	a_{\sigma}\cdot \dfrac{2u_{2\sigma}^i}{u_{\lambda}}=\frac{a_{\sigma}^3u_{2\sigma}^{i-1}}{a_{\lambda}}
	\end{gather}
		Note: If $i=0$ in the first equation or $j=0$ in the 2nd-4th equations then we get a negative exponent in a denominator. When this happens that means the product is $0$.
	\begin{proof}
		After clearing denominators (multiplying by the denominators in the right hand side), all but the third equation reduce to the Gold relation. The third instead reduces to \eqref{xs1}.
	\end{proof}
	
	\item The following relations compute the product of $u_{2\sigma}$ with the other generators:
	\begin{gather}
		u_{2\sigma}\cdot \dfrac{2s}{a_{\sigma}a_{\lambda}^iu_{\lambda}^j}=\frac{a_{\sigma}s}{a_{\lambda}^{i+1}u_{\lambda}^{j-1}}\\
	\label{Second}	u_{2\sigma}\cdot \dfrac{2a_{\lambda}^i}{a_{\sigma}}=a_{\sigma}a_{\lambda}^{i-1}u_{\lambda}
	\end{gather}

\begin{proof}
	Both reduce to the Gold as usual. For the second equation, multiplication by $a_{\sigma}$ is an isomorphism as can be seen directly from the right hand side.
\end{proof}
		
	\item The following relations compute the product of $a_{\lambda}$ with the other generators:
		\begin{gather}
	a_{\lambda}\cdot \dfrac{2s}{a_{\sigma}u_{2\sigma}^iu_{\lambda}^j}=\frac{x_{1,1}}{u_{2\sigma}^iu_{\lambda}^j}\\
\label{alambdauu}a_{\lambda}\cdot \dfrac{2u_{2\sigma}^i}{u_{\lambda}}=a_{\sigma}^2u_{2\sigma}^{i-1}
	\end{gather}	
	\begin{proof}
		Both reduce to the Gold after clearing denominators.
	\end{proof}

\item The following relations compute the product of $u_{\lambda}$ with the other generators:
	 \begin{gather}
	u_{\lambda}\cdot \dfrac{2s}{a_{\sigma}u_{2\sigma}^ia_{\lambda}^j}=0\\
	u_{\lambda}\cdot \dfrac{w_3}{a_{\sigma}^iu_{2\sigma}^j}=0\\
	u_{\lambda}\cdot 	\dfrac{a_{\sigma}^i}{a_{\lambda}^j}=0\textup{ , }i\ge 3
	\end{gather}
	
	\begin{proof}The first two relations are deduced as follows:
		\begin{gather}
	u_{\lambda}\cdot \dfrac{2s}{a_{\sigma}u_{2\sigma}^ia_{\lambda}^j}=0\cdot \frac{x_{1,1}}{a_{\sigma}^2u_{2\sigma}^{i-1}a_{\lambda}^{j}}\iff 
u_{\lambda}a_{\sigma}\dfrac{2s}{u_{2\sigma}}=0\impliedby 2a_{\sigma}=0\\
u_{\lambda}\cdot \dfrac{w_3}{a_{\sigma}^iu_{2\sigma}^j}=0\cdot \frac{a_{\lambda}w_3}{a_{\sigma}^{i+2}u_{2\sigma}^{j-1}}\iff 
a_{\sigma}^2u_{\lambda}\dfrac{w_3}{u_{2\sigma}}=0\impliedby Gold
		\end{gather}
For the last relation, if $j\ge 2$ the homology group in the degree of the product is $0$ so we may assume $j=1$ and further that $i=3$ (we can factor higher powers of $a_{\sigma}$ out of the quotient). Then,
		\begin{gather}
		u_{\lambda}\cdot 	\dfrac{a_{\sigma}^3}{a_{\lambda}}=0\cdot a_{\sigma}u_{2\sigma}\iff a_{\sigma}^3u_{\lambda}=0\impliedby Gold
		\end{gather}
	\end{proof}
	
	\item The remaining relations involving $s$ are:
	\begin{gather}
	\dfrac{s}{u_{2\sigma}^iu_{\lambda}^j}\cdot \dfrac{2a_{\lambda}}{a_{\sigma}}=\frac{x_{1,1}}{u_{2\sigma}^iu_{\lambda}^j}\\
	\dfrac{2s}{a_{\sigma}u_{2\sigma}^ia_{\lambda}^j}\cdot \dfrac{2a_{\lambda}}{a_{\sigma}}=0\\
	\dfrac{2s}{a_{\sigma}u_{2\sigma}^ia_{\lambda}^ju_{\lambda}^k}\cdot \dfrac{2a_{\lambda}}{a_{\sigma}}=2\cdot \frac{s}{u_{2\sigma}^{i+1}a_{\lambda}^ju_{\lambda}^{k-1}}
	\end{gather}
	\begin{proof}The first reduces to \eqref{xs2} as usual. The second reduces to the relation we just proved, while the final one reduces to \eqref{Second}:
		\begin{gather}
		\dfrac{2s}{a_{\sigma}u_{2\sigma}^ia_{\lambda}^j}\cdot \dfrac{2a_{\lambda}}{a_{\sigma}}=0\cdot \frac{x_{1,1}}{a_{\sigma}u_{2\sigma}^ia_{\lambda}^j}\iff 2s\frac{2a_{\lambda}}{a_{\sigma}}=0\\
		\dfrac{2s}{a_{\sigma}u_{2\sigma}^ia_{\lambda}^ju_{\lambda}^k}\cdot \dfrac{2a_{\lambda}}{a_{\sigma}}=2\cdot \frac{s}{u_{2\sigma}^{i+1}a_{\lambda}^ju_{\lambda}^{k-1}}\iff \frac{2s}{a_{\sigma}u_{\lambda}}u_{2\sigma}\frac{2a_{\lambda}}{a_{\sigma}}=2s\impliedby\eqref{Second}
		\end{gather}
	\end{proof}

	\item The remaining relations involving $x_{1,1}$ are:
	\begin{equation}
	\begin{gathered}
	\frac{x_{1,1}}{a_{\sigma}u_{2\sigma}^ia_{\lambda}^j}\cdot \frac{2u_{2\sigma}}{u_{\lambda}}=2\frac{s}{u_{2\sigma}^ia_{\lambda}^{j}}\\
	\frac{x_{1,1}}{a_{\sigma}^iu_{2\sigma}^ja_{\lambda}^k}\cdot \frac{2u_{2\sigma}}{u_{\lambda}}=\frac{x_{1,1}}{a_{\sigma}^{i-2}u_{2\sigma}^ja_{\lambda}^{k+1}}	\end{gathered}\label{x112uu}
	\end{equation}

	\begin{proof}		For the first relation we perform a denominator exchange:
		\begin{equation}
		\begin{gathered}
		\frac{x_{1,1}}{a_{\sigma}u_{2\sigma}^ia_{\lambda}^j}\cdot \frac{2u_{2\sigma}}{u_{\lambda}}=2\frac{s}{u_{2\sigma}^ia_{\lambda}^{j}}\iff %\frac{x_{1,1}}{a_{\sigma}} \frac{2u_{2\sigma}}{u_{\lambda}}=2s	\iff 
		\frac{x_{1,1}}{a_{\sigma}u_{\lambda}} 2u_{2\sigma}=2s\iff \frac{s}{u_{2\sigma}}2u_{2\sigma}=s\\
		\frac{x_{1,1}}{a_{\sigma}^iu_{2\sigma}^ja_{\lambda}^k}\cdot \frac{2u_{2\sigma}}{u_{\lambda}}=\frac{x_{1,1}}{a_{\sigma}^{i-2}u_{2\sigma}^ja_{\lambda}^{k+1}}\iff \frac{x_{1,1}}{a_{\sigma}^2} a_{\lambda} \frac{2u_{2\sigma}}{u_{\lambda}}=x_{1,1}\impliedby\eqref{alambdauu}
		\end{gathered}
		\end{equation}
	\end{proof}

	\item The remaining relations involving $w_3$ are:
	\begin{gather}
	\dfrac{w_3}{a_{\sigma}^iu_{2\sigma}^j}\cdot \dfrac{a_{\sigma}^3}{a_{\lambda}^k}=\frac{x_{1,1}}{a_{\sigma}^{i-1}u_{2\sigma}^{j}a_{\lambda}^{k-1}}\\
		\frac{w_3}{a_{\sigma}^iu_{2\sigma}^j}\cdot \frac{2u_{2\sigma}}{u_{\lambda}}=\frac{x_{1,1}}{a_{\sigma}^iu_{2\sigma}^j}
	\end{gather}
	
	\begin{proof}For the first relation we may assume $i\ge 1$ (otherwise we have a transfer) and then the equality is implied by \eqref{xw}. The second relation is implied by \eqref{Eq1}.
	\end{proof}
	\item The remaining relation involving $u_{\lambda}/u_{2\sigma}$ is:
	\begin{gather}\label{uuaa}
	\frac{u_{\lambda}}{u_{2\sigma}}\frac{a_{\sigma}^3}{a_{\lambda}}=0
	\end{gather}
	\begin{proof}
		Follows immediately from the Gold.
	\end{proof}
	\item The remaining relations involving $(2a_{\lambda})/a_{\sigma}$ are:
	\begin{gather}\label{aaaa}
	\dfrac{2a_{\lambda}}{a_{\sigma}}\cdot \dfrac{2a_{\lambda}}{a_{\sigma}}=2\frac{a_{\lambda}u_{\lambda}}{u_{2\sigma}}\\
		\dfrac{2a_{\lambda}}{a_{\sigma}}\cdot \dfrac{a_{\sigma}^3}{a_{\lambda}}=0\\
		\dfrac{2a_{\lambda}^i}{a_{\sigma}}\cdot \dfrac{2u_{2\sigma}^j}{u_{\lambda}}=0
	\end{gather}
	
	\begin{proof}For the first two:
		\begin{gather}
		\dfrac{2a_{\lambda}}{a_{\sigma}}\cdot \dfrac{2a_{\lambda}}{a_{\sigma}}=2\frac{a_{\lambda}u_{\lambda}}{u_{2\sigma}}\iff \dfrac{2a_{\lambda}}{a_{\sigma}}\cdot (u_{2\sigma}\dfrac{2a_{\lambda}}{a_{\sigma}})=2a_{\lambda}u_{\lambda}\impliedby \eqref{Second}\\
		\dfrac{2a_{\lambda}}{a_{\sigma}}\cdot \dfrac{a_{\sigma}^3}{a_{\lambda}}=0\cdot a_{\sigma}^2\impliedby 2a_{\lambda}\cdot  \dfrac{a_{\sigma}^3}{a_{\lambda}}=0\cdot a_{\sigma}^3
		\end{gather}
For the third we may assume $j=1$ and then
		\begin{gather}\label{2aa2uu}
		\dfrac{2a_{\lambda}^i}{a_{\sigma}}\cdot \dfrac{2u_{2\sigma}}{u_{\lambda}}=0\cdot a_{\sigma}a_{\lambda}^{i-1}\iff 2a_{\lambda}^i\cdot \dfrac{2u_{2\sigma}}{u_{\lambda}}=0\cdot a_{\sigma}^2a_{\lambda}^{i-1}\impliedby \eqref{alambdauu}
		\end{gather}
	\end{proof}
	
	\item The remaining relations involving $2u_{2\sigma}/u_{\lambda}$ are:	\begin{gather}\label{u2u2ulul}
	\dfrac{2u_{2\sigma}^i}{u_{\lambda}}\cdot \dfrac{2u_{2\sigma}^j}{u_{\lambda}}=\frac{4u_{2\sigma}^{i+j}}{u_{\lambda}^2}+\frac{a_{\sigma}^4u_{2\sigma}^{i+j-2}}{a_{\lambda}^2}\\
		\frac{a_{\sigma}^i}{a_{\lambda}^k}\cdot \frac{2u_{2\sigma}^j}{u_{\lambda}}=\frac{a_{\sigma}^{i+2}u_{2\sigma}^{j-1}}{a_{\lambda}^{k+1}}
	\end{gather}
	\begin{proof}We explained how to get the first relation at the end of subsection \ref{Subtle}. For the second, note that the element in the left-hand side is $2$-torsion so no torsion-free generator can appear in the right-hand side. Thus we have
				\begin{gather}
		\frac{a_{\sigma}^i}{a_{\lambda}^k}\cdot \frac{2u_{2\sigma}^j}{u_{\lambda}}=\frac{a_{\sigma}^{i+2}u_{2\sigma}^{j-1}}{a_{\lambda}^{k+1}}\iff a_{\sigma}^ia_{\lambda} \frac{2u_{2\sigma}^j}{u_{\lambda}}=a_{\sigma}^{i+2}u_{2\sigma}^{j-1}\impliedby \eqref{alambdauu}
		\end{gather}
	\end{proof}
	
\end{itemize}

\phantom{1}\smallbreak

\begin{small}
	\noindent  \textsc{Department of Mathematics, University of Chicago}\\
	\textit{E-mail:} \verb|nickg@math.uchicago.edu|\\
	\textit{Website:} \href{http:://math.uchicago.edu/~nickg}{math.uchicago.edu/$\sim$nickg}
\end{small}

\end{document}